\documentclass[12pt,twoside,a4paper]{article}
\title{Higher dimensional Zariski decompositions}
\author{S\'ebastien Boucksom}

\usepackage{latexsym}
\begin{document}
\maketitle

\newcommand{\er}{\mathbf{R}}
\newcommand{\ku}{\mathbf{Q}}
\newcommand{\pet}{\mathbf{P}}
\newcommand{\ze}{\mathbf{Z}}
\newcommand{\co}{\mathbf{C}}
\newcommand{\oh}{\mathcal{O}}
\newcommand{\ih}{\mathcal{I}}
\newcommand{\peff}{\mathcal{E}}
\newcommand{\ka}{\mathcal{K}}
\newcommand{\bigc}{\mathcal{B}}
\newcommand{\bka}{\mathcal{BK}}
\newcommand{\mka}{\mathcal{MK}}
\newcommand{\nef}{\mathcal{N}}
\newcommand{\mnef}{\mathcal{MN}}
\newcommand{\pos}{\mathcal{P}}
\newcommand{\cone}{\peff^{\star}}
\newcommand{\ach}{H^{1,1}(X,\er)}
\newcommand{\ddbar}{\partial\overline{\partial}}
\newcommand{\dach}{H^{1,1}_{\ddbar}(X,\er)}
\newcommand{\nq}{NS(X)\otimes\ku}
\newcommand{\nr}{NS(X)_{\er}}
\newcommand{\ep}{\varepsilon}
\newcommand{\de}{\delta}
\newcommand{\ti}{\widetilde}
\newtheorem{theo}{Theorem}[section]
\newtheorem{conj}{Conjecture}
\newtheorem{cor}[theo]{Corollary}
\newtheorem{defi}[theo]{Definition}
\newtheorem{prop}[theo]{Proposition}
\newtheorem{lem}[theo]{Lemma}

$Author's$ $address$: Institut Fourier, 100 rue des Maths, BP74, 38402 Saint-Martin d'H\`eres Cedex, France.\\
$e-mail$: sbouckso@ujf-grenoble.fr\\
$Abstract$: using currents with minimal singularities, we construct minimal multiplicities for a real pseudo-effective $(1,1)$-class $\alpha$ on a compact complex $n$-fold $X$, which are the local obstructions to the numerical effectivity of $\alpha$. The negative part of $\alpha$ is then defined as the real effective divisor $N(\alpha)$ whose multiplicity along a prime divisor $D$ is just the generic multiplicity of $\alpha$ along $D$, and we get in that way a divisorial Zariski decomposition of $\alpha$ into the sum of a class $Z(\alpha)$ which is nef in codimension 1 and the class of its negative part $N(\alpha)$, which is exceptional in the sense that it is very rigidly embedded in $X$. The positive parts $Z(\alpha)$ generate a modified nef cone, and the pseudo-effective cone is shown to be locally polyhedral away from the modified nef cone, with extremal rays generated by exceptional divisors. We then treat the case of a surface and a hyper-K\"ahler manifold in some detail: under the intersection form (resp. the Beauville-Bogomolov form), we characterize the modified nef cone and the exceptional divisors; our divisorial Zariski decomposition is orthogonal, and is thus a rational decomposition, which fact accounts for the usual existence statement of a Zariski decomposition on a projective surface, which is thus extended to the hyper-K\"ahler case. Finally, we explain how the divisorial Zariski decomposition of (the first Chern class of) a big line bundle on a projective manifold can be characterized in terms of the asymptotics of the linear series $|kL|$ as $k\to\infty$. 

\section{Introduction}
It is known since the pioneering work of O.Zariski [Zar62] that the study of the linear series $|kL|$ where $L$ is a line bundle on a projective surface can be reduced to the case where $L$ is numerically effective (nef). The more precise result obtained by Zariski is that any effective $\ku$-divisor $D$ on a projective surface $X$ can be uniquely decomposed into a sum $D=P+N$ where $P$ is a nef $\ku$-divisor (the positive part), $N=\sum a_j D_j$ is an effective $\ku$-divisor (the negative part) such that the Gram matrix $(D_i\cdot D_j)$ is negative definite, and $P$ is orthogonal to $N$ with respect to the intersection form. Zariski shows that the natural inclusion $H^0(kP)\to H^0(kL)$ is necessarily an isomorphism in that case, relating the decomposition to the original problem.\\
The proof of the uniqueness in this decomposition shows that the negative part $N$ only depends on the class $\{D\}$ of $D$ in the N\'eron-Severi group $NS(X)$, so that $\{D\}\mapsto\{P\}$ yields a map from part of the pseudo-effective cone to the nef cone, which we want to study geometrically.\\
Building upon the construction by J.-P.Demailly of metrics with minimal singularities on a pseudo-effective line bundle $L$ over a compact complex $n$-fold, we define the minimal multiplicity $\nu(\alpha,x)$ of an arbitrary real pseudo-effective $(1,1)$-class $\alpha$ on a compact complex $n$-fold $X$ at some point $x\in X$. This multiplicity $\nu(\alpha,x)$ is the local obstruction at $x$ to the numerical effectivity of $\alpha$, and we then get the negative part of such a class $\alpha$ by setting $N(\alpha)=\sum\nu(\alpha,D)D$, where $D$ ranges over the prime divisors of $X$ and $\nu(\alpha,D)$ is the generic multiplicity of $\alpha$ along $D$ (cf. section 3). This negative part $N(\alpha)$ is an effective $\er$-divisor which is exceptional in the sense that it is very rigidly imbedded in $X$. When $X$ is a surface, the divisors we obtain in that way are exactly the effective $\er$-divisors whose support $D_1,...,D_r$ have negative definite Gram matrix $(D_i\cdot D_j)$.\\
The difference $Z(\alpha):=\alpha-\{N(\alpha)\}$ is a real $(1,1)$-class on $X$ which we call the Zariski projection of $\alpha$. It is not a nef class, but is somehow nef in codimension 1. More precisely, we define the modified nef cone of a K\"ahler $n$-fold to be the closed convex cone generated by the classes in $\ach$ which can be written as the push-forward of a K\"ahler class by a modification. We then show that the Zariski projection $Z(\alpha)$ of a pseudo-effective class $\alpha$ belongs to this modified nef cone. The decomposition $\alpha=Z(\alpha)+\{N(\alpha)\}$ we call the divisorial Zariski decomposition, and it is just induced by the Siu decomposition of a positive current with minimal singularities in $\alpha$ when the latter is big. For such a big class, we give a criterion to recognize a decomposition $\alpha=p+\{N\}$ into a modified nef and big class and the class of an effective real divisor as the divisorial Zariski decomposition of $\alpha$, in terms of the non-K\"ahler locus of $p$ (cf. section 3.5)\\
The geometric picture is now as follows: the pseudo-effective cone of a compact complex $n$-fold $X$ is locally polyhedral away from the modified nef cone, with extremal rays that write $\er_+\{D\}$ for some exceptional prime $D$ of $X$. The Zariski projection $Z$ yields a projection from the pseudo-effective cone to the modified nef cone parallel to these exceptional rays, which map is concave (in some sense) and homogeneous, but not continuous up to the boundary of the pseudo-effective cone in general. The fibre $Z^{-1}(p)$ of $Z$ above a modified nef class $p$ is a countable union of simplicial cones generated by exceptional families of primes.\\
When $X$ is a surface, a modified nef class is just a nef class; when $\alpha$ is the class of an effective $\ku$-divisor $D$ on a projective surface, the divisorial Zariski decomposition of $\alpha$ is just the original Zariski decomposition of $D$. More generally, we show that the divisorial decomposition of a pseudo-effective class $\alpha$ on a K\"ahler surface is the unique orthogonal decomposition of $\alpha$ into the sum of a modified nef class and the class of an exceptional (in some sense) effective $\er$-divisor. This fact accounts for the rationality of the Zariski decomposition on a surface, meaning that the negative part $N$ is rational when $D$ is.\\
An interesting fact is that much of the well-known case of a surface carries on to the case where $X$ is a compact hyper-K\"ahler manifold. Using the quadratic Beauville-Bogomolov form on $\ach$ and deep results due to D.Huybrechts, we can prove the following facts: a family of primes is exceptional in our sense iff their Gram matrix is negative definite. In particular, a prime is exceptional iff it has negative square, and this forces it to be uniruled. The modified nef cone of a hyper-K\"ahler manifold is just the dual cone to the pseudoeffective cone, which is also the closure of the so-called birational (or bimeromorphic) K\"ahler cone. Finally, the divisorial Zariski decomposition is the unique orthogonal decomposition into a modified nef class and an exceptional divisor. In particular, the divisorial Zariski decomposition is also rational in that case.\\ 
In a last part, we explain how to tackle the above constructions in a more algebraic fashion. When $L$ is a big divisor on a projective manifold, we prove that the divisorial Zariski projection of $L$ is the only decomposition $L=P+N$ into real divisors with $P$ modified nef and $H^0(\lfloor kP\rfloor)=H^0(kL)$ for every $k$. The minimal multiplicities of $\{L\}$ (and thus its negative part) can be recovered from the asymptotic behaviour of the sections of $kL$. The case of a general pseudo-effective line bundle $L$ is then handled by approximating it by $L+\ep A$, where $A$ is ample.\\
The methods used in this paper are mostly ``transcendental'', since we heavily rely on the theory of currents, but the results we aim at definitely belong to algebraic geometry, and we have thus tried to make the paper legible for more algebraically inclined readers by providing in the first section a rather detailed account of the tools we need afterwards.               

\section{Technical preliminaries}
\subsection{$\ddbar$-cohomology}
When $X$ is an arbitrary complex manifold, the $\ddbar$-lemma of
K\"ahler geometry does not hold, and it is thus better to work with
$\ddbar$-cohomology. We will just need the $(1,1)$-cohomology space
$H^{1,1}_{\ddbar}(X,\co)$, which is defined as the quotient of the space of $d$-closed smooth $(1,1)$-forms modulo the $\ddbar$-exact ones. The real structure on the space of forms induces a real structure on $H^{1,1}_{\ddbar}(X,\co)$, and we denote by $\dach$ the space of real points.\\

The canonical map from $H^{1,1}_{\ddbar}(X,\co)$ to the quotient of the space of $d$-closed $(1,1)$-currents modulo the $\ddbar$-exact ones is injective (because, for any degree $0$ current $f$, $\ddbar f$ is smooth iff $f$ is), and is also surjective: given a closed $(1,1)$-current $T$, one can find a locally finite open covering $U_j$ of $X$ such that $T=\ddbar f_j$ is $\ddbar$-exact on $U_j$. If $\rho_j$ is a partition of unity associated to $U_j$ and $f:=\sum\rho_jf_j$, then $T-\ddbar f$ is smooth. Indeed, on $U_i$, it is just $\ddbar(\sum_j\rho_j(f_i-f_j)$, and each $f_i-f_j$ is smooth since it is even pluri-harmonic. As a consequence, a class $\alpha\in H^{1,1}_{\ddbar}(X,\co)$ can be seen as an affine space of closed $(1,1)$-currents. We denote by $\{T\}\in H^{1,1}_{\ddbar}(X,\co)$ the class of the current $T$. Remark that $i\ddbar$ is a real operator (on forms or currents), so that if $T$ is a real closed $(1,1)$-current, its class $\{T\}$ lies in $\dach$ and consists in all the closed currents $T+i\ddbar\varphi$ where $\varphi$ is a real current of degree $0$.\\

When $X$ is furthermore compact, it can be shown that $H^{1,1}_{\ddbar}(X,\co)$ is finite dimensional. The operator $\ddbar$ from smooth functions to smooth closed (1,1)-forms is thus an operator between Fr\'echet spaces with finite codimensional range; it therefore has closed range, and the quotient map $\theta\mapsto\{\theta\}$ from smooth closed $(1,1)$-forms to $H^{1,1}_{\ddbar}(X,\co)$ endowed with its unique finite-dimensional complex vector space Hausdorff topology is thus continuous and open.

\subsection{General facts about currents}
\subsubsection{Siu decomposition}
Let $T$ be a closed positive current of bidegree $(p,p)$ on a complex
$n$-fold $X$. We denote by $\nu(T,x)$ its Lelong number at a point
$x\in X$. The Lelong super-level sets are defined by $E_c(T):=\{x\in
X,\nu(T,x)\geq c\}$ for $c>0$, and a well known result of Y.T.Siu
[Siu74] asserts that $E_c(T)$ is an analytic subset of $X$, of codimension at least $p$. As a consequence, for any analytic subset $Y$ of $X$, the generic Lelong number of $T$ along $Y$, defined by 
$$\nu(T,Y):=\inf\{\nu(T,x),x\in Y\},$$
is also equal to $\nu(T,x)$ for a very general $x\in Y$. It is also true that, for any irreducible analytic subset $Y$ of codimension $p$ in $X$, the current \\$T-\nu(T,Y)[Y]$ is positive. The symbol $[Y]$ denotes the integration current on $Y$, which is defined by integrating test forms on the smooth locus of $Y$. Since $E_+(T):=\cup_{c>0}E_c(T)$ is a countable union of $p$-codimensional analytic subsets, it contains an at most countable family $Y_k$ of $p$-codimensional irreducible analytic subsets. By what we have said, $T-\nu(T,Y_1)[Y_1]-...-\nu(T,Y_k)[Y_k]$ is a positive current for all $k$, thus the series $\sum_{k\geq 0}\nu(T,Y_k)[Y_k]$ converges, and we have 
$$T=R+\sum_k\nu(T,Y_k)[Y_k]$$ 
for some closed positive $(p,p)$-current $R$ such that each $E_c(R)$ has codimension $>p$. The decomposition above is called the Siu decomposition of the closed positive $(p,p)$-current $T$. Since $\nu(T,Y)=0$ if $Y$ is a $p$-codimensional subvariety not contained in $E_+(T)$, it makes sense to write $\sum_k\nu(T,Y_k)[Y_k]=\sum\nu(T,Y)[Y]$, where the sum is implicitely extended over all $p$-codimensional irreducible analytic subsets $Y\subset X$.\\ 

\subsubsection{Almost positive currents}
A real $(1,1)$-current $T$ on a complex manifold $X$ is said to be almost positive if $T\geq\gamma$ holds for some smooth real $(1,1)$-form $\gamma$. Let $T\geq\gamma$ be a closed almost positive $(1,1)$-current. On a small enough open set $U$ with coordinates $z=(z_1,...,z_n)$, we write $T=\ddbar\varphi$ where $\varphi$ is a degree $0$ current. Since $\gamma+Ci\ddbar|z|^2$ is a positive $(1,1)$-form on $U$ for $C>0$ big enough, we get that $i\ddbar(\varphi+C|z|^2)$ is positive, which means that $\varphi+C|z|^2$ is (the current associated to) a (unique) pluri-subharmonic function on $U$. A locally integrable function $\varphi$ on $X$ such that $i\ddbar\varphi$ is almost positive is called an almost pluri-subharmonic function, and is thus locally equal to a pluri-subharmonic function modulo a smooth function.\\

The Lelong number $\nu(T,x)$ of a closed almost positive $(1,1)$-current $T$ can be defined as $\nu(T+Ci\ddbar|z|^2,x)$ as above, since this does not depend on the smooth function $C|z|^2$. Consequently, the Siu decomposition of $T$ can also be constructed, and writes $T=R+\sum\nu(T,D)[D]$, where $D$ ranges over the prime divisors of $X$, and $R$ is a closed almost positive $(1,1)$-current. In fact, we have $R\geq\gamma$ as soon as $T\geq\gamma$ for a smooth form $\gamma$.\\

\subsubsection{Pull-back of a current}
When $f:Y\to X$ is a $surjective$ holomorphic map between compact complex manifolds and $T$ is a closed almost positive $(1,1)$-current on $X$, it is possible to define its pull back $f^{\star}T$ by $f$ using the analogue of local equations for divisors: write $T=\theta+i\ddbar\varphi$ for some smooth form $\theta\in\{T\}$. $\varphi$ is then an almost pluri-subharmonic function, thus locally a pluri-subharmonic function modulo $\mathcal{C}^{\infty}$. One defines $f^{\star}T$ to be $f^{\star}\theta+i\ddbar f^{\star}\varphi$, as this is easily seen to be independent of the choices made. Of course, we then have $\{f^{\star}T\}=f^{\star}\{T\}$.

\subsubsection{Gauduchon metrics and compactness}
On any compact complex $n$-fold $X$, there exists a Hermitian metric
$\omega$ such that $\omega^{n-1}$ is $\ddbar$-closed. Such a metric is
called a Gauduchon metric. As a consequence, for every smooth real $(1,1)$-form $\gamma$, the quotient map $T\mapsto\{T\}$ from the set of closed $(1,1)$-currents $T$ with $T\geq\gamma$ to $\dach$ is proper. Indeed, the mass of the positive current $T-\gamma$ is controled by $\int(T-\gamma)\wedge\omega^{n-1}$, and $\int T\wedge\omega^{n-1}=\{T\}\cdot\{\omega\}$ only depends on the class of $T$. The result follows by the weak compactness of positive currents with bounded mass. Another consequence is that the kernel of $T\mapsto\{T\}$ meets the cone of closed positive $(1,1)$-currents at the origin only.\\ 

\subsubsection{Cycles as currents}
One can associate to any effective $p$-codimensional $\er$-cycle $Y=a_1Y_1+...+a_rY_r$ a closed positive $(p,p)$-current $[Y]=a_1[Y_1]+...+a_r[Y_r]$, called the integration current on $Y$. The map $Y\mapsto[Y]$ so defined is injective, and a result of Thie says that the Lelong number $\nu([Y],x)$ is just the multiplicity of $Y$ at $x$. Consequently, we shall drop the brackets in $[Y]$ when no confusion is to be feared, and write for instance $T=R+\sum\nu(T,D)D$ for a Siu decomposition, because this is more in the spirit of this work. 

\subsection{Cones in the $\ddbar$-cohomology}
We now assume that $X$ is compact, and fix some reference Hermitian
 form $\omega$ (i.e. a smooth positive definite $(1,1)$-form). A
 cohomology class \\$\alpha\in\dach$ is said to be pseudo-effective
 iff it contains a positive current;\\
$\alpha$ is nef (numerically effective) iff, for each $\ep>0$,
 $\alpha$ contains a smooth form $\theta_{\ep}$ with
 $\theta_{\ep}\geq-\ep\omega$;\\
$\alpha$ is big iff it contains a K\"ahler current, i.e. a closed
 $(1,1)$-current $T$ such that $T\geq\ep\omega$ for $\ep>0$ small
 enough. Finally, $\alpha$ is a K\"ahler class iff it contains a
 K\"ahler form (note that a smooth K\"ahler current is the same thing
 as a K\"ahler form).\\
Since any two Hermitian forms $\omega_1$, $\omega_2$ are commensurable
 (i.e. $C^{-1}\omega_2\leq\omega_1\leq C\omega_2$ for some $C>0$),
 these definitions do not depend on the choice of $\omega$.\\
 
The set of pseudo-effective classes is a closed convex cone $\peff\subset\dach$, called the pseudo-effective cone. It has compact base, because so is the case of the cone of closed positive $(1,1)$-currents. Similarly, one defines the nef cone $\nef$ (a closed convex cone), the big cone $\bigc$ (an open convex cone), and the K\"ahler cone $\ka$ (an open convex cone). We obviously have the inclusions 
$$\ka\subset\bigc\subset\peff$$
and
$$\ka\subset\nef\subset\peff.$$
By definition, $X$ is a K\"ahler manifold iff its K\"ahler cone $\ka$
is non-empty. Similarly (but this is a theorem, cf. [DP01]) $X$ is a Fujiki manifold (i.e. bimeromorphic to a K\"ahler manifold) iff its big cone $\bigc$ is non-empty (see also the proof of proposition 2.3 below). If $X$ is K\"ahler, $\ka$ is trivially the interior $\nef^0$ of the nef cone. Similarly, if $X$ is Fujiki, $\bigc$ is trivially the interior $\peff^0$ of the pseudo-effective cone.\\
We will now and then denote by $\geq$ the partial order relation on $\dach$ induced by the convex cone $\peff$.

\subsection{The N\'eron-Severi space}
Given a line bundle $L$ on $X$, each smooth Hermitian metric $h$ on
$L$ locally writes as $h(x,v)=|v|^2e^{-2\varphi(x)}$ for some smooth
local weight $\varphi$; the curvature form
$\Theta_h(L):=\frac{i}{\pi}\ddbar\varphi$ is a globally defined real
$(1,1)$-form, whose class in $\dach$ we denote by $c_1(L)$, the first
Chern class of $L$. We write $dd^c=\frac{i}{\pi}\ddbar$ for short. A
singular Hermitian metric $h$ on $L$ is by definition a metric
$h=h_{\infty}e^{-2\varphi}$, where $h_{\infty}$ is a smooth Hermitian
metric on $L$ and the weight $\varphi$ is a locally integrable
function. The curvature current of $h$ is defined as
$\Theta_h(L):=\Theta_{h_{\infty}}(L)+dd^c\varphi$; it also lies in
$c_1(L)$. Conversely, given a smooth Hermitian metric $h_{\infty}$ on
$L$, any closed real $(1,1)$-current $T$ in $c_1(L)$ can be written
(by definition) as $T=\Theta_{h_{\infty}}(L)+dd^c\varphi$. But
$\varphi$ is just a degree $0$ current $a$ $priori$. However,
$\varphi$ is automatically $L^1_{loc}$ in case $T$ is almost positive (cf. section 2.2.2), thus each almost positive current $T$ in $c_1(L)$
is the curvature form of a singular Hermitian metric on $L$.\\

The image of the homomorphism Pic$(X)\to\dach$ $L\mapsto c_1(L)$ is called the N\'eron-Severi group, denoted by $NS(X)$. It is a free $\ze$-module, whose rank is denoted by $\rho(X)$, and called the Picard number of $X$. The real N\'eron-Severi space $\nr$ is just the real subspace of dimension $\rho(X)$ in $\dach$ generated by $NS(X)$. Kodaira's embedding theorem can be formulated as follows: $X$ is a projective manifold iff the intersection of the K\"ahler cone $\ka$ with $\nr$ is non-empty. Similarly, $X$ is a Moishezon manifold (i.e. bimeromorphic to a projective manifold) iff the intersection of the big cone $\bigc$ with $\nr$ is non-empty (cf. [DP01]).

\subsection{Currents with analytic singularities}
\subsubsection{Definition}
A closed almost positive $(1,1)$-current $T$ on a compact complex $n$-fold $X$ is said to have analytic singularities (along a subscheme $V(\ih)$ defined by a coherent ideal sheaf $\ih$) if there exists some $c>0$ such that $T$ is locally congruent to $\frac{c}{2}dd^c\log(|f_1|^2+...+|f_k|^2)$ modulo smooth forms, where $f_1,...,f_k$ are local generators of $\ih$. $T$ is thus smooth outside the support of $V(\ih)$, and it is an immediate consequence of the Lelong-Poincar\'e formula that $\sum\nu(T,D)D$ is just $c$ times the divisor part of the scheme $V(\ih)$. If we first blow-up $X$ along $V(\ih)$ and then resolve the singularities, we get a modification $\mu:\ti{X}\to X$, where $\ti{X}$ is a compact complex manifold, such that $\mu^{-1}\ih$ is just $\oh(-D)$ for some effective divisor $D$ upstairs. The pull-back $\mu^{\star}T$ clearly has analytic singularities along $V(\mu^{-1}\ih)=D$, thus its Siu decomposition writes
$$\mu^{\star}T=\theta+cD$$
where $\theta$ is a smooth $(1,1)$-form. If $T\geq\gamma$ for some smooth form $\gamma$, then $\mu^{\star}T\geq\mu^{\star}\gamma$, and thus $\theta\geq\mu^{\star}\gamma$. This operation we call a resolution of the singularities of $T$. 

\subsubsection{Regularization(s) of currents}
We will need two basic types of regularizations (inside a fixed cohomology class) for closed $(1,1)$-currents, both due to J.-P.Demailly. 
\begin{theo}[Dem82, Dem92]
Let $T$ be a closed almost positive $(1,1)$-current on a compact
complex manifold $X$, and fix a Hermitian form $\omega$. Suppose that
$T\geq\gamma$ for some smooth real $(1,1)$-form $\gamma$ on $X$. Then:\\ 

(i) There exists a sequence of smooth forms $\theta_k$ in $\{T\}$ which converges weakly to $T$, and such that $\theta_k\geq\gamma-C\lambda_k\omega$ where $C>0$ is a constant depending on the curvature of $(T_X,\omega)$ only, and $\lambda_k$ is a decreasing sequence of continuous functions such that $\lambda_k(x)\to\nu(T,x)$ for every $x\in X$.\\

(ii) There exists a sequence $T_k$ of currents with analytic singularities in $\{T\}$ which converges weakly to $T$, such that $T_k\geq\gamma-\ep_k\omega$ for some sequence $\ep_k>0$ decreasing to $0$, and such that $\nu(T_k,x)$ increases to $\nu(T,x)$ uniformly with respect to $x\in X$.  

\end{theo}

Point (ii) enables us in particular to approximate a K\"ahler current $T$ inside its cohomology class by K\"ahler currents $T_k$ with analytic singularities, with a very good control of the singularities. A big class therefore contains plenty of K\"ahler currents with analytic singularities.       

\subsection{Intersection of currents}
Just as cycles, currents can be intersected provided their singular
sets are in an acceptable mutual position. Specifically, let $T$ be a
closed positive $(1,1)$-current on a complex manifold $X$. Locally, we
have $T=dd^c\varphi$ with $\varphi$ a pluri-subharmonic function,
which is well defined modulo a pluri-harmonic (hence smooth)
function. We therefore get a globally well-defined unbounded locus
$L(T)$, which is the complement of the open set of points near which
$\varphi$ is locally bounded. Assume now that $T_1$, $T_2$ are two
closed positive $(1,1)$-currents such that $L(T_j)$ is contained in an
analytic set $A_j$ (which may be $X$); locally, we write
$T_j=dd^c\varphi_j$ with $\varphi_j$ a pluri-subharmonic function. If
$A_1\cap A_2$ has codimension at least $2$, then it is shown in
[Dem92] that $\varphi_1dd^c\varphi_2$ has locally finite mass, and
that $dd^c\varphi_1\wedge dd^c\varphi_2:=dd^c(\varphi_1
dd^c\varphi_2)$ yields a globally defined closed positive
$(2,2)$-current, denoted by $T_1\wedge T_2$. It is also true that
$T_1\wedge T_2$ lies in the product cohomology class
$\{T_1\}\cdot\{T_2\}\in H^{2,2}_{\ddbar}(X,\er)$.\\
We will only need the following two special cases: if $T_1$ is a closed positive $(1,1)$-current with analytic singularities along a subscheme of codimension at least $2$, then $T_1\wedge T_2$ exists for every closed positive $(1,1)$-current $T_2$. \\
If $D_1$ and $D_2$ are two distinct prime divisors, then $[D_1]\wedge
[D_2]$ is a well defined closed positive $(2,2)$-current. Since its
support is clearly contained in the set-theoretic intersection
$D_1\cap D_2$ (whose codimension is at least $2$), we have
$[D_1]\wedge [D_2]=\sum a_j[Y_j]$, where the $Y_j$'s are the
components of $D_1\cap D_2$. In fact, it can be shown that $\sum
a_jY_j$ is just the $2$-cycle associated to the scheme-theoretic
intersection $D_1\cap D_2$, thus $[D_1]\wedge[D_2]$ is just the integration current associated to the cycle $D_1\cdot D_2$.  

\subsection{The modified nef cone}
For our purposes, we need to introduce a new cone in $\dach$, which is somehow the cone of classes that are nef in codimension 1. Let $X$ be a compact complex $n$-fold, and $\omega$ be some reference Hermitian form. 

\begin{defi}[Modified nef and K\"ahler classes]Let $\alpha$ be a class in $\dach$. 

\noindent(i) $\alpha$ is said to be a modified K\"ahler class iff it contains a K\"ahler current $T$ with $\nu(T,D)=0$ for all prime divisors $D$ in $X$.

\noindent(ii) $\alpha$ is said to be a modified nef class iff, for every $\ep>0$, there exists a closed $(1,1)$-current $T_{\ep}$ in $\alpha$ with $T_{\ep}\geq-\ep\omega$ and $\nu(T_{\ep},D)=0$ for every prime $D$. 
\end{defi}
This is again independent of the choice of $\omega$ by commensurability of the Hermitian forms. The set of modified K\"ahler classes is an open convex cone called the modified K\"ahler cone and denoted by $\mka$. Similarly, we get a closed convex cone $\mnef$, the modified nef cone. Using the Siu decomposition, we immediately see that $\mka$ is non-empty iff the big cone $\bigc$ is non-empty, in which case $\mka$ is just the interior of $\mnef$.\\ 
{\bf Remark 1}: upon regularizing the currents using (ii) of theorem 2.1, we can always assume that the currents involved in the definition have analytic singularities along a subcheme of codimension at least 2.\\
{\bf Remark 2}: the modified nef cone of a compact complex surface is just its nef cone (cf. section 4.2.1).\\
{\bf Remark 3}: just as for nef classes, one cannot simply take
$\ep=0$ in the definition of a modified nef class. We recall the
example given in [DPS94]: there exists a ruled surface $X$ over an
elliptic curve such that $X$ contains an irreducible curve $C$ with the following property: the class $\{C\}\in\dach$ is nef, but contains only one positive current, which is of course the integration current $[C]$. 

\noindent The following proposition gives a more ``algebraic'' characterization of $\mka$, which also explains the (seemingly dumb) terminology.

\begin{prop} A class $\alpha$ lies in $\mka$ iff there exists a modification\\
$\mu:\ti{X}\to X$ and a K\"ahler class $\ti{\alpha}$ on $\ti{X}$ such that
$\alpha=\mu_{\star}\ti{\alpha}$.
\end{prop}
$Proof$: the argument is adapted from [DP01], theorem 3.4. If $\ti{\omega}$ is a K\"ahler form on $\ti{X}$ and $\omega$ is our reference Hermitian form on $X$, then $\mu^{\star}\omega\leq C\ti{\omega}$ for some $C>0$, since $\ti{X}$ is compact. Since $\mu$ is a modification, we have $\mu_{\star}\mu^{\star}\omega=\omega$, so we get $T:=\mu_{\star}\ti{\omega}\geq C^{-1}\omega$, and $T$ is thus a K\"ahler current. Since the singular values of $\mu$ are in codimension at least 2, we immediately see that $\nu(T,D)=0$ for every prime divisor $D$ in $X$, and $\{T\}=\mu_{\star}\{\omega\}$ lies in $\mka$ as desired. Conversely, if $\alpha\in\mka$ is represented by a K\"ahler current $T$ with $\nu(T,D)=0$ for all $D$, there exists by (ii) of theorem 2.1 a K\"ahler current $T_k$ in $\alpha$ with analytic singularities along a subscheme $V_k$ with $\nu(T_k,D)\leq\nu(T,D)$, so that $V_k$ has no divisor component. We select a resolution of the singularities of $T_k$ $\mu:\ti{X}\to X$, and write $\mu^{\star}T_k=\theta+F$, where $\theta$ is a smooth form and $F$ is an effective $\er$-divisor. Since $T_k\geq\ep\omega$ for $\ep>0$ small enough, we get that $\theta\geq\mu^{\star}\ep\omega$. Denoting by $E_1,...,E_r$ the $\mu$-exceptional prime divisors on $\ti{X}$, it is shown in [DP01], lemma 3.5, that one can find $\de_1,...,\de_r>0$ small enough and a closed smooth $(1,1)$-form $\tau$ in $\{\de_1E_1+...+\de_rE_r\}$ such that $\mu^{\star}\ep\omega-\tau$ is positive definite everywhere. It follows that $\theta-\tau$ is a K\"ahler form upstairs. Now, we have 
$$\alpha=\mu_{\star}\{T_k\}=\mu_{\star}\{\theta-(\de_1E_1+...+\de_rE_r)\}=\mu_{\star}\{\theta-\tau\},$$
since $E_j$ is $\mu$-exceptional and so is $F$ because $\mu_{\star}F$ is an effective divisor contained in the scheme $V_k$; this concludes the proof of proposition 2.3.

\noindent That a modified nef classe is somehow nef in codimension 1 is reflected in the following
\begin{prop} If $\alpha$ is a modified K\"ahler (resp. nef) class, then $\alpha_{|D}$ is big (resp. pseudo-effective) for every prime divisor $D\subset X$. 
\end{prop}
$Proof$: if $\alpha$ is a modified nef class and $\ep>0$ is given,
choose a current $T_{\ep}\geq-\ep\omega$ in $\alpha$ with analytic
singularities in codimension at least 2. Locally, we have
$\omega\leq Cdd^c|z|^2$ for some $C>0$, thus
$T_{\ep}+\ep Cdd^c|z|^2$ writes as $dd^c\varphi_{\ep}$, where
$\varphi_{\ep}$ is pluri-subharmonic and is not identically $-\infty$
on $D$. Thus the restriction $(\varphi_{\ep})_{|D}$ is
pluri-subharmonic, and $(T_{\ep}+\ep Cdd^c|z|^2)_{|D}$ is a well
defined closed positive current. It follows that $(T_{\ep})_{|D}$ is a
well defined almost positive current on $D$, with
$(T_{\ep})_{|D}\geq-\ep C\omega_{|D}$. This certainly implies that $\alpha_{|D}$ is pseudo-effective. The case $\alpha\in\mka$ is treated similarly. 

\subsection{Currents with minimal singularities}
Let $\varphi_1$, $\varphi_2$ be two almost pluri-subharmonic functions on a compact complex manifold $X$. Then, following [DPS00], we say that $\varphi_1$ is less singular than
$\varphi_2$ (and write $\varphi_1\preceq\varphi_2$) if we have
$\varphi_2\leq\varphi_1+C$ for some constant $C$. We denote by $\varphi_1\approx\varphi_2$ the equivalence relation generated by the pre-order relation $\preceq$. Note that $\varphi_1\approx\varphi_2$ exactly means that $\varphi_1=\varphi_2$ mod $L^{\infty}$.\\ 
When $T_1$ and $T_2$ are two closed almost positive $(1,1)$-currents on $X$, we can also compare their singularities in the following fashion: write $T_i=\theta_i+dd^c\varphi_i$ for
$\theta_i\in\{T_j\}$ a smooth form and $\varphi_i$ an almost pluri-subharmonic function. Since any $L^1_{loc}$ function $f$ with $dd^cf$ smooth is itself smooth, it is easy to check that $\varphi_i$ does not depend on the choices made up to equivalence of singularities, and we compare the singularities of the $T_i$'s by comparing those of the $\varphi_i$'s.\\
Let now $\alpha$ be a class in $\dach$ and $\gamma$ be a smooth real $(1,1)$-form, and denote by $\alpha[\gamma]$ the set of closed almost positive $(1,1)$-currents $T$ lying in $\alpha$ with $T\geq\gamma$. It is a (weakly) compact and convex subset of the space of $(1,1)$-currents. We endow it with the pre-order relation $\preceq$ defined above. For any family $T_j$, $j\in J$ of elements of $\alpha[\gamma]$, we claim that there exists an infimum $T=\inf_{j\in J}T_j$ in $(\alpha[\gamma],\preceq)$, which is therefore unique up to equivalence of singularities. The proof is pretty straightforward: fix a smooth form $\theta$ in $\alpha$, and write $T_j=\theta+dd^c\varphi_j$ for some quasi pluri-subharmonic functions $\varphi_j$. Since $X$ is compact, $\varphi_j$ is bounded from above; therefore, upon changing $\varphi_j$ into $\varphi_j-C_j$, we may assume that $\varphi_j\leq 0$ for all $ j\in J$. We then take $\varphi$ to be the upper semi-continuous upper enveloppe of the $\varphi_j$'s, $j\in J$, and set $T:=\theta+dd^c\varphi$. It is immediate to check that $T\preceq T_j$ for all $ j$, and that for every $S\in\alpha[\gamma]$, $S\preceq T_j$ for all $ j$ implies that $S\preceq T$. We should maybe explain why $T\geq\gamma$: locally, we can choose coordinates $z=(z_1,...,z_n)$ and a form $q(z)=\sum\lambda_j|z_j|^2$ such that $dd^cq\leq\gamma$ and $dd^cq$ is arbitrarily close to $\gamma$. Writing $\theta=dd^c\psi$ for some smooth local potential $\psi$, the condition $\theta+dd^c\varphi_j\geq\gamma$ implies that $\psi+\varphi_j-q$ is pluri-subharmonic. The upper enveloppe $\psi+\varphi-q$ is thus also pluri-subharmonic, which means that $T=\theta+dd^c\varphi\geq dd^cq$; letting $dd^cq$ tend to $\gamma$, we get $T\geq\gamma$, as desired.\\

Since any two closed almost positive currents with equivalent singularities have the same Lelong numbers, the Lelong numbers of $\inf T_j$ do not depend on the specific choice of the current. In fact, it is immediate to check from the definitions that
$$\nu(\inf_{j\in J} T_j,x)=\inf_{j\in J}\nu(T_j,x).$$ \\

As a particular case of the above construction, there exists a closed
almost positive $(1,1)$-current $T_{\min,\gamma}\in\alpha[\gamma]$ which is a least element in $(\alpha[\gamma],\preceq)$. $T_{\min,\gamma}$ is well defined modulo $dd^cL^{\infty}$, and we call it a current with minimal singularities in $\alpha$, for the given lower bound $\gamma$. 
When $\gamma=0$ and $\alpha$ is pseudo-effective, we just write $T_{\min}=T_{\min,0}$, and call it a positive current with minimal singularities in $\alpha$. It must be noticed that, even for a big class $\alpha$, $T_{\min}$ will be a K\"ahler current only in the trivial case:
\begin{prop}A pseudo-effective class $\alpha$ contains a positive current with
minimal singularities $T_{\min}$ which is a K\"ahler current iff $\alpha$ is a K\"ahler class.
\end{prop}
$Proof$: we can write $T_{\min}=\theta+dd^c\varphi$ with $\theta$ a smooth
form. If $T_{\min}$ is K\"ahler, then so is
$\ep\theta+(1-\ep)T_{\min}=\theta+dd^c(1-\ep)\varphi$ for $\ep>0$ small enough.
We therefore get $\varphi\preceq(1-\ep)\varphi$ by minimality, that is:
$(1-\ep)\varphi\leq\varphi+C$ for some constant $C$. But this shows that
$\varphi$ is bounded, and thus $T_{\min}$ is a K\"ahler current with
identically zero Lelong numbers. Using (i) of theorem 2.1, we can therefore regularize it into a K\"ahler form inside its cohomology class, qed.\\

Finally, we remark that a positive current with minimal singularities in a pseudo-effective class is generally non-unique (as a current), as the example of a K\"ahler class already shows.

\section{The divisorial Zariski decomposition}
In this section $X$ denotes a compact complex $n$-fold, and $\omega$ is a reference Hermitian form, unless otherwise specified.

\subsection{Minimal multiplicities and non-nef locus}
When $\alpha\in\dach$ is a pseudo-effective class, we want to
introduce minimal multiplicities $\nu(\alpha,x)$, which measure the
obstruction to the numerical effectivity of $\alpha$. For each
$\ep>0$, let $T_{\min,\ep}=T_{\min,\ep}(\alpha)$ be a current with
minimal singularities in $\alpha[-\ep\omega]$ (cf. section 2.8 for the
notation). We then introduce the following
\begin{defi}[Minimal multiplicities] The minimal multiplicity at $x\in X$ of the pseudo-effective class $\alpha\in\dach$ is defined as
$$\nu(\alpha,x):=\sup_{\ep>0}\nu(T_{\min,\ep},x).$$
\end{defi}
The commensurability of any two Hermitian forms shows that the definition does not depend on $\omega$.
\noindent When $D$ is a prime divisor, we define the generic minimal multiplicity of $\alpha$ along $D$ as 
$$\nu(\alpha,D):=\inf\{\nu(\alpha,x),x\in D\}.$$
We then have $\nu(\alpha,D)=\sup_{\ep>0}\nu(T_{\min,\ep},D)$, and $\nu(\alpha,D)=\nu(\alpha,x)$ for the very general $x\in D$.
\begin{prop} Let $\alpha\in\dach$ be a pseudo-effective class. 

\noindent(i) $\alpha$ is nef iff $\nu(\alpha,x)=0$ for every $x\in X$.

\noindent(ii) $\alpha$ is modified nef iff $\nu(\alpha,D)=0$ for every prime $D$.
\end{prop}
$Proof$: if $\alpha$ is nef (resp. modified nef), $\alpha[-\ep\omega]$
contains by definition a smooth form (resp. a current $T_{\ep}$ with
$\nu(T_{\ep},D)=0$ for every prime $D$). We thus have
$\nu(T_{\min,\ep},x)=0$ (resp. $\nu(T_{\min,\ep},D)=0$) for every
$\ep>0$, and thus $\nu(\alpha,x)=0$
(resp. $\nu(\alpha,D)=0$). Conversely, if $\nu(\alpha,x)=0$ for every
$x\in X$, applying (i) of theorem 2.1 to $T_{\min,\ep}$, we see that
$\nu(T_{\min,\ep},x)=0$ for every $x\in X$ implies that
$\alpha[-\ep'\omega]$ contains a smooth form for every $\ep'>\ep$, and
$\alpha$ is thus nef. Finally, if $\nu(\alpha,D)=0$ for every prime
$D$, we have $\nu(T_{\min,\ep},D)=0$ for every prime $D$. Since
$T_{\min,\ep}$ lies in $\alpha[-\ep\omega]$, $\alpha$ is modified nef
by the very definition.\\

\noindent In view of proposition 3.2, we propose the 
\begin{defi}[Non-nef locus] The non-nef locus of a pseudo-effective class $\alpha\in\dach$ is defined by 
$$E_{nn}(\alpha):=\{x\in X,\nu(\alpha,x)>0\}.$$
\end{defi}
Recall that the set $E_+(T):=\{x\in X,\nu(T,x)>0\}$ is a countable union of closed analytic subsets for every closed almost positive $(1,1)$-current $T$. Since $E_{nn}(\alpha)=\cup_{\ep>0}E_+(T_{\min,\ep})$, the non-nef locus is also a countable union of closed analytic subsets. We do not claim however that each super-level set $\{x\in X,\nu(\alpha,x)\geq c\}$ $(c>0)$ is an analytic subset (this is most certainly not true in general). Using results of M.Paun, proposition 3.2 generalizes as follows:
\begin{prop} A pseudo-effective class $\alpha$ is nef iff $\alpha_{|Y}$ is nef for every irreducible analytic subset $Y\subset E_{nn}(\alpha)$.
\end{prop}
$Proof$: since the restriction of a nef class to any analytic subset is nef, one direction is clear. To prove the converse, we cannot directly apply the results of M.Paun [Pau98], since we allow slightly negative currents, so we sketch the proof to show that it immediately carries on to our situation. We quote without proof the following results: 
\begin{lem}[Gluing lemma] Let $\alpha$ be any class in $\dach$, and $Y_1$, $Y_2$ two analytic subsets of $X$. If $\alpha_{|Y_i}$ is nef $(i=1,2)$, then $\alpha_{|Y_1\cup Y_2}$ is nef.
\end{lem}

\begin{lem}[Extension lemma] Let $\theta$ be a closed smooth $(1,1)$-form and $\gamma$ be any smooth $(1,1)$-form. Let $Y\subset X$ be an analytic subset, and assume that 
$$\theta_{|Y}+dd^c\varphi\geq\gamma_{|Y}$$
for some smooth function $\varphi$ on $Y$. Then, for every $\ep>0$, there exists a neighbourhood $V$ of $Y$ and a smooth function $\varphi_{\ep}$ on $V$ such that 
$$\theta+dd^c\varphi_{\ep}\geq\gamma-\ep\omega$$
on $V$. 
\end{lem}
We now select once for all a smooth form $\theta$ in $\alpha$. By applying (ii) of theorem 2.1 to $T_{\min,\ep}\geq-\ep\omega$, we can select a closed $(1,1)$-current with analytic singularities $T^{(1)}_{\ep}=\theta+dd^c\varphi^{(1)}_{\ep}$ such that $T^{(1)}_{\ep}\geq-2\ep\omega$ and $\nu(T^{(1)}_{\ep},x)\leq\nu(T_{\min,\ep},x)$. Let $Y_{\ep}$ be the analytic subset along which $T^{(1)}_{\ep}$ is singular; we have $Y_{\ep}\subset E_+(T_{\min,\ep})\subset E_{nn}(\alpha)$. Since $\alpha$ is nef by assumption on every component of $Y_{\ep}$, using the above two lemma, we can find a neighbourhood $V_{\ep}$ of $Y_{\ep}$ and a smooth function $\varphi^{(2)}_{\ep}$ on $V_{\ep}$ such that $\theta+dd^c\varphi^{(2)}_{\ep}\geq-\ep\omega$ on $V_{\ep}$. We choose a smaller neighbourhood $W_{\ep}$ of $Y_{\ep}$ with $\overline{W_{\ep}}\subset V_{\ep}$, and we then set 
$$\varphi^{(3)}_{\ep}:=\cases{\varphi^{(1)}_{\ep}& on $X-W_{\ep}$,\cr
\max_{\eta}(\varphi_{\ep}^{(2)}-C_{\ep},\varphi^{(1)}_{\ep})& on $\overline{W_{\ep}}$\cr}$$ 
where $\max_{\eta}(x,y):=\max\star\rho_{\eta}$ denotes a regularized
maximum function obtained by convolution with a regularizing kernel
$\rho_{\eta}$ ($\eta$ is chosen so small that $\max_{\eta}(x,y)=x$
when $y<x-1/2$), and $C_{\ep}$ is a positive constant, large enough to
achieve $\varphi^{(1)}_{\ep}\geq\varphi_{\ep}^{(2)}-C_{\ep}+1$ near
$\partial W_{\ep}$ (we use that $\varphi^{(1)}_{\ep}$ is smooth away
from $Y_{\ep}$, hence locally bounded near $\partial W_{\ep}$). The
two parts to be glued then coincide near $\partial W_{\ep}$, thus
$\varphi^{(3)}_{\ep}$ is smooth. Since both
$\theta+dd^c\varphi^{(1)}_{\ep}$ and
$\theta+dd^c\varphi_{\ep}^{(2)}-C_{\ep}$ are greater than
$-2\ep\omega$, the gluing property of pluri-subharmonic functions yields that $\theta+dd^c\varphi^{(3)}_{\ep}\geq-2\ep\omega$. Since this is true for every $\ep>0$, this shows that $\alpha$ is indeed nef.   

We now investigate the continuity of $\alpha\mapsto\nu(\alpha,x)$ and $\nu(\alpha,D)$:
\begin{prop} For every $x\in X$ and every prime $D$, the maps
  $\peff\to\er$ $\alpha\mapsto\nu(\alpha,x)$ and $\nu(\alpha,D)$ are
  convex, homogeneous. They are continuous on the interior $\peff^0$,
  and lower semi-continuous on the whole of $\peff$.
\end{prop}
$Proof$: let $\alpha$, $\beta$ be two pseudo-effective classes. If $T_{\min,\ep}(\alpha)$ and $T_{\min,\ep}(\beta)$ are currents with minimal singularities in $\alpha[-\ep\omega]$ and $\beta[-\ep\omega]$ respectively, then \\$T_{\min,\ep}(\alpha)+T_{\min,\ep}(\alpha)$ belongs to $(\alpha+\beta)[-2\ep\omega]$, thus 
$$\nu(T_{\min,2\ep}(\alpha+\beta),x)\leq\nu(T_{\min,\ep}(\alpha),x)+\nu(T_{\min,\ep}(\beta),x)\leq\nu(\alpha,x)+\nu(\beta,x).$$
We infer from this $\nu(\alpha+\beta,x)\leq\nu(\alpha,x)+\nu(\beta,x)$, and a similar sub-additivity property for $\nu(\cdot,D)$ is obtained along the same lines. Since the homogeneity of our two maps is obvious, the convexity also follows.\\
The quotient map $\theta\mapsto\{\theta\}$ from the Fr\'echet space of
closed smooth real $(1,1)$-forms to $\dach$ is surjective, thus
open. If $\alpha_k\in\dach$ is a sequence of pseudo-effective classes
converging to $\alpha$ and $\ep>0$ is given, we can thus find a smooth
form $\theta_k\in\alpha-\alpha_k$ for each $k$ big enough such that
\\$-\ep\omega\leq\theta_k\leq\ep\omega$. The current
$T_{\min,\ep}(\alpha_k)+\theta_k$ then lies in $\alpha[-2\ep\omega]$,
and thus
$\nu(T_{\min,2\ep}(\alpha),x)\leq\nu(T_{\min,\ep}(\alpha_k),x)\leq\nu(\alpha_k,x)$,
for each $k$ big enough. We infer from this that
$\nu(T_{\min,2\ep}(\alpha),x)\leq\liminf_{k\to\infty}\nu(\alpha_k,x)$
for each $\ep>0$, hence
$\nu(\alpha,x)\leq\liminf_{k\to\infty}\nu(\alpha_k,x)$, by taking the
supremum of the left hand-side for $\ep>0$. This means that
$\alpha\mapsto\nu(\alpha,x)$ is lower semi-continuous, and similarly for $\nu(\alpha,D)$, just replacing $x$ by $D$ in the above proof.\\
Finally, the restrictions of our maps to $\peff^0$ are continuous as any convex map on an open convex subset of a finite dimensional vector space is.
   
\begin{prop}Let $\alpha\in\dach$ be a pseudo-effective class, and
  $T_{\min}$ be a positive current with minimal singularities in
  $\alpha$.

(i) We always have $\nu(\alpha,x)\leq\nu(T_{\min},x)$ and
$\nu(\alpha,D)\leq\nu(T_{\min},D)$.

(ii) When $\alpha$ is furthermore big, we have $\nu(\alpha,x)=\nu(T_{\min},x)$ and
$\nu(\alpha,D)=\nu(T_{\min},D)$.

\end{prop}
$Proof$: since $T_{\min}$ belongs to $\alpha[-\ep\omega]$ for every
$\ep>0$, $\nu(\alpha,x)\leq\nu(T_{\min},x)$ follows for every $x\in
X$, for any pseudo-effective class $\alpha$. If $\alpha$ is
furthermore big, we can choose a K\"ahler current $T$ in $\alpha$ with
$T\geq\omega$ for some Hermitian form $\omega$. If $T_{\min,\ep}$ is a
current with minimal singularities in $\alpha[-\ep\omega]$, then
$(1-\ep)T_{\min,\ep}+\ep T$ is a positive current in $\alpha$, and
thus $\nu((1-\ep)T_{\min,\ep}+\ep T,x)\geq\nu(T_{\min}
,x)$ by
  minimality of $T_{\min}$, from which we infer
$$(1-\ep)\nu(\alpha,x)+\ep\nu(T,x)\geq\nu(T_{\min},x).$$
We thus get the converse inequality $\nu(\alpha,x)\geq\nu(T_{\min},x)$
  by letting $\ep\to 0$. The case of $\nu(\alpha,D)$ is similar.

\subsection{Definition of the divisorial Zariski decomposition}
Let $\alpha\in\dach$ be again a pseudo-effective class, and choose a
positive current with minimal singularities $T_{\min}$ in
$\alpha$. Since $\nu(\alpha,D)\leq\nu(T_{\min},D)$ for every prime $D$
by proposition 3.8, the series of currents $\sum\nu(\alpha,D)[D]$ is convergent, since it is dominated by $\sum\nu(T_{\min},D)[D]$. 
\begin{defi}[Divisorial Zariski decomposition] The negative part of a pseudo-effective class $\alpha\in\dach$ is defined as $N(\alpha):=\sum\nu(\alpha,D)[D]$. The Zariski projection of $\alpha$ is $Z(\alpha):=\alpha-\{N(\alpha)\}$. We call the decomposition $\alpha=Z(\alpha)+\{N(\alpha)\}$ the divisorial Zariski decomposition of $\alpha$.
\end{defi}
It is certainly highly desirable that the negative part $N(\alpha)$ of a pseudo-effective class be a divisor, i.e. that $\nu(\alpha,D)=0$ for almost every prime $D$. We will see in section 3.3 that it is indeed the case. For the time being, we concentrate on the Zariski projection, which we see as a map $Z:\peff\to\peff$. 
\begin{prop}Let $\alpha\in\dach$ be a pseudo-effective class. Then:

\noindent(i) Its Zariski projection $Z(\alpha)$ is a modified nef class.

\noindent(ii) We have $Z(\alpha)=\alpha$ iff $\alpha$ is modified nef. 

\noindent(iii) $Z(\alpha)$ is big iff $\alpha$ is.

\noindent(iv) If $\alpha$ is not modified nef, then $Z(\alpha)$
belongs to the boundary $\partial\mnef$ of the modified nef cone.
\end{prop}
$Proof$:(i) Let $T_{\min,\ep}$ be as before a current with minimal
singularities in $\alpha[-\ep\omega]$, and consider its Siu
decomposition
$T_{\min,\ep}=R_{\ep}+\sum\nu(T_{\min,\ep},D)[D]$. First, we claim
that $N_{\ep}:=\sum\nu(T_{\min,\ep},D)[D]$ converges weakly to
$N(\alpha)$ as $\ep$ goes to $0$. For any smooth form $\theta$ of
bidimension $(1,1)$, $\theta+C\omega^{n-1}$ is a positive form for
$C>0$ big enough. Every such $\theta$ is thus the difference of two
positive forms, and it is enough to show that $\int
N_{\ep}\wedge\theta\to\int N(\alpha)\wedge\theta$ for every smooth
positive form $\theta$. But $\int
N_{\ep}\wedge\theta=\sum\nu(T_{\min,\ep},D)\int[D]\wedge\theta$ is a
convergent series whose general term
$\nu(T_{\min,\ep},D)\int[D]\wedge\theta$ converges to
$\nu(\alpha,D)\int[D]\wedge\theta$ as $\ep\to 0$ and is dominated by
$\nu(T_{\min},D)\int[D]\wedge\theta$; since
$\sum\nu(T_{\min},D)\int[D]\wedge\theta\leq\int T_{\min}\wedge\theta$
converges, our claim follows by dominated convergence.\\
In particular, the class $\{N_{\ep}-N(\alpha)\}$ converges to
zero. Since the map $\theta\mapsto\{\theta\}$ is open on the space of
smooth closed $(1,1)$-form, we can find a sequence
$\theta_k\geq-\de_k\omega$ of smooth forms with
$\theta_k\in\{N_{\ep_k}-N(\alpha)\}$ for some sequences $\ep_k<<\de_k$
going to zero. It remains to notice that $T_k:=R_{\ep_k}+\theta_k$ is
a current in $Z(\alpha)$ with $T_k\geq-(\ep_k+\de_k)\omega$ and
$\nu(T_k,D)=0$ for every prime $D$. Since $\ep_k+\de_k$ converges to
zero, $Z(\alpha)$ is modified nef by definition.   

\noindent(ii) Since $N(\alpha)=\sum\nu(\alpha,D)[D]$ is a closed positive $(1,1)$-current, it is zero iff its class $\{N(\alpha)\}\in\dach$ is. The assertion is thus just a reformulation of (ii) in proposition 3.2.

\noindent(iii) If $Z(\alpha)$ is big, then of course $\alpha=Z(\alpha)+\{N(\alpha)\}$ is also big, as the sum of a big class and a pseudo-effective one. If conversely $\alpha$ is big, it contains a K\"ahler current $T$, whose Siu decomposition we write $T=R+\sum\nu(T,D)[D]$. Note that $R$ is a K\"ahler current since $T$ is; since $T$ belongs to $\alpha[-\ep\omega]$ for every $\ep>0$, we have $\nu(T,D)\geq\nu(\alpha,D)$, and $R+\sum(\nu(T,D)-\nu(\alpha,D))[D]$ is thus a K\"ahler current in $Z(\alpha)$ as desired.

\noindent(iv) Assume that $Z(\alpha)$ belongs to the interior
$\mnef^0$ of the modified nef cone. By proposition 3.2, we have to see
that $\nu(\alpha,D)=0$ for every prime $D$. Suppose therefore that
$\nu(\alpha,D_0)>0$ for some prime $D_0$. The class
$Z(\alpha)+\ep\{D_0\}$ has to lie in the open cone $\mnef^0$ for $\ep$
small enough, thus we can write for $0<\ep<\nu(\alpha,D_0)$:
$$\alpha=(Z(\alpha)+\ep\{D_0\})+(\nu(\alpha,D_0)-\ep)\{D_0\}+\{\sum_{D\neq
  D_0}\nu(\alpha,D)D\}.$$
We deduce that
$\nu(\alpha,D_0)\leq\nu(Z(\alpha)+\ep\{D_0\},D_0)+(\nu(\alpha,D_0)-\ep)$.
Indeed, the class $\{D_0\}$ (resp. $\{\sum_{D\neq D_0}\nu(\alpha,D)D\}$)
has minimal multiplicity $\leq 1$ (resp. 0) along $D_0$, because so is
the
generic Lelong numbers of the positive current $[D_0]$
(resp. $\sum_{D\neq D_0}\nu(\alpha,D)[D]$) along $D_0$. Now, we also
have $\nu(Z(\alpha)+\ep\{D_0\},D_0)=0$ since $Z(\alpha)+\ep\{D_0\}$ is modified nef by assumption, hence the
contradiction $\nu(\alpha,D_0)\leq\nu(\alpha,D_0)-\ep$. 

\begin{prop}(i) The map $\alpha\mapsto N(\alpha)$ is convex and
  homogeneous on $\peff$. It is continuous on the interior of the pseudo-effective cone.

(ii) The Zariski projection $Z:\peff\to\mnef$ is concave and homogeneous. It is continuous on the interior of $\peff$. 
\end{prop}
$Proof$: we have already noticed that $\nu(\alpha+\beta,D)\leq\nu(\alpha,D)+\nu(\beta,D)$ for every prime $D$ and every two pseudo-effective classes $\alpha$, $\beta$. This implies that $N(\alpha+\beta)\leq N(\alpha)+N(\beta)$. Homogeneity is obvious, and the first assertion follows. To show continuity, it is enough as above to show that $\alpha\mapsto\int N(\alpha)\wedge\theta$ is continuous on $\peff^0$ for every positive form $\theta$. But the latter map is convex, and thus continuous on $\peff^0$ as any convex map on an open convex subset of a finite dimensional vector space is. (ii) is now an obvious consequence of (i) and the relation $Z(\alpha)=\alpha-\{N(\alpha)\}$. 
\subsection{Negative part and exceptional divisors} 
If $A=D_1,...,D_r$ is a finite family of prime divisors, we denote by $V_+(A)\subset\dach$ the closed convex cone generated by the classes $\{D_1\},...,\{D_r\}$. Every element of $V_+(A)$ writes $\alpha=\{E\}$ for some effective $\er$-divisor supported by the $D_j$'s. Since $[E]$ is a positive current in $\alpha$,  we have $N(\alpha)\leq E$, and thus $Z(\alpha)$ can be represented by the effective $\er$-divisor $E-N(\alpha)$, which is also supported by the $D_j$'s. We conclude: $V_+(A)$ is stable under the Zariski projection $Z$. In particular, we have $Z(V_+(A))=0$ iff $V_+(A)$ meets $\mnef$ at $0$ only. 
\begin{defi}[Exceptional divisors](i) A family $D_1,...,D_q$ of prime divisors is said to be an exceptional family iff the convex cone generated by their cohomology classes meets the modified nef cone $\mnef$ at $0$ only.

\noindent(ii) An effective $\er$-divisor $E$ is said to be exceptional iff its prime components constitute an exceptional family. 
\end{defi}

We have the following 
\begin{prop}(i) An effective $\er$-divisor $E$ is exceptional iff $Z(\{E\})=0$.

\noindent(ii) If $E$ is an exceptional effective $\er$-divisor, we have $E=N(\{E\})$. 

\noindent(iii) If $D_1,...,D_q$ is an exceptional family of primes, then their classes $\{D_1\},...,\{D_q\}$ are linearly independent in $\nr\subset\dach$. In particular, the length of the exceptional families of primes is uniformly bounded by the Picard number $\rho(X)$.
\end{prop}
$Proof$: (i) let $A=D_1,...,D_r$ denote the family of primes
supporting $E$, and choose a Gauduchon metric $\omega$ (cf. section
2.2.4). Since $\omega^{n-1}$ is $\ddbar$-closed, $\int
Z(\alpha)\wedge\omega^{n-1}$ is well defined, and defines a map
$\peff\to\er$ $\alpha\mapsto\int Z(\alpha)\wedge\omega^{n-1}$, which
is concave and homogeneous (by proposition 3.11), and everywhere
non-negative. The restriction of this map to $V_+(A)$ shares the same
properties, and the class $\alpha:=\{E\}$ is a point in the relative
interior of the convex cone $V_+(A)$ at which $\int
Z(\alpha)\wedge\omega^{n-1}=0$. By concavity, we thus get $\int
Z(\alpha)\wedge\omega^{n-1}=0$ for every $\alpha\in V_+(A)$, and thus
$Z(\alpha)=0$ for every such $\alpha\in V_+(A)$, qed.\\
(ii) When $E$ is exceptional, we have both $E\geq N(\{E\})$ (because
the positive current $[E]$ lies in the class $\{E\}$) and $\{E\}=\{N(\{E\})\}$ (because $Z(\{E\})=0$). Since a closed positive current which yields zero in $\dach$ is itself zero, we get the result.\\
(iii) Since $D_1,...,D_q$ are linearly independent in Div$(X)\otimes\er$, the assertion is equivalent to the fact that the quotient map $D\mapsto\{D\}$ is injective on the $\er$-vector space of divisors generated by the $D_j$'s. But this is easy: if $E=\sum a_j D_j$ lies in the kernel, we can write $E=E_+-E_-$ with $E_+$ and $E_-$ effective such that $\{E_+\}=\{E_-\}$. By (iii), we get $E_+=E_-$, whence $E=0$, qed. 

We state as a theorem the following important consequences of (iii):
\begin{theo}(i) For every pseudo-effective class $\alpha\in\peff$, the
  negative part $N(\alpha)$ is an exceptional effective $\er$-divisor supported by at most $\rho(X)$ primes.

\noindent(ii) $X$ carries at most countably many exceptional primes.

\noindent(iii) The exceptional fiber $Z^{-1}(0)$ is contained in $\nr$, and is a union of at most countably many simplicial cones over exceptional families of primes. 
\end{theo}
$Proof$: (i) We have $Z(\alpha)\geq Z(Z(\alpha))+Z(\{N(\alpha)\})$,
and $Z(Z(\alpha))=Z(\alpha)$ by proposition 3.10, thus
$Z(\{N(\alpha)\})=0$. We immediately deduce from this that any family
of primes $D_1,...,D_r$ such that $\nu(\alpha,D_j)>0$ for every $j$ is
an exceptional family, and the assertion follows from (iii) of
proposition 3.13.\\ 
(ii) We just have to notice that $D\mapsto\{D\}$ is injective on the set of exceptional primes, and maps into the lattice $NS(X)\subset\nr$.\\
(iii) Since $\{A\}$ is a linearly independent set for every exceptional family of primes $A$, we see that $V_+(A)=\sum_{D\in A}\er_+\{D\}$ is a simplicial cone. It remains to observe that $\alpha$ lies in the exceptional fiber $Z^{-1}(0)$ iff $\alpha=\{N(\alpha)\}$, thus $Z^{-1}(0)$ is covered by the simplicial cones $V_+(A)$.\\

We will see in section 4.3 that a family $D_1,...,D_q$ of primes on a
surface is exceptional iff the Gram matrix $(D_i\cdot D_j)$ is
negative definite, i.e. iff $D_1,...,D_q$ can all be blown down to
points by a modification towards an analytic surface (singular in
general). On a general compact complex $n$-fold $X$, an exceptional
divisor is still very rigidly embedded in $X$:
\begin{prop} If $E$ is an exceptional effective $\er$-divisor, then
  its class $\{E\}$ contains but one positive current, which is
  $[E]$. In particular, when $E$ is rational, its Kodaira-Iitaka
  dimension $\kappa(X,E)$ is zero.
\end{prop}
$Proof$: if $T$ is a positive current in $\{E\}$, we have
$\nu(T,D)\geq\nu(\{E\},D)$ for every prime $D$. Using the Siu
decomposition of $T$, we thus see that
$T\geq\sum\nu(\{E\},D)D=N(\{E\})=E$, since $E$ is exceptional. But we
also have $\{T\}=\{E\}$, hence $T=E$, as was to be shown. To get the
last point, let $D$ be an element of the linear system $|kE|$ for some
integer $k>0$ such that $kE$ is Cartier. The positive current
$\frac{1}{k}[D]$ then lies in $\{E\}$, thus we have $[D]=k[E]$ as
currents, hence $D=kE$ as divisors. This shows that $h^0(kE)=1$ for
each $k>0$, qed.  

\subsection{Discontinuities of the Zariski projection}
It is remarkable that the Zariski projection $Z$ is not continuous in general up to the boundary $\partial\peff$. 

\begin{prop} If $X$ carries infinitely many exceptional primes, then
  the Zariski projection $Z:\peff\to\mnef$ is not continuous.
\end{prop}
$Proof$: we use the following
\begin{lem} If $D_k$ is an infinite sequence of divisors, the rays $\er_+\{D_k\}\subset\peff$ can accumulate on $\mnef$ only.
\end{lem}  
$Proof$: suppose that $t_k\{D_k\}$ converges to some non-zero
$\alpha\in\peff$ (for $t_k>0$). For each prime $D$, we then have
$D_k\neq D$ and thus $\nu(t_k\{D_k\},D)=t_k\nu(\{D_k\},D)=0$ for
infinitely many $k$, because the family $D_k$ is infinite. By lower
semi-continuity (proposition 3.7) we deduce $\nu(\alpha,D)=0$ for every
prime $D$, i.e. $\alpha$ is modified nef (by proposition 3.2).   

\noindent Assume now that an infinite sequence of exceptional prime
divisors $D_k$ exists. Since $\peff$ has compact base, upon extracting
a subsequence, we can assume that $t_k\{D_k\}$ converges to some
non-zero $\alpha\in\peff$ (with $t_k>0$ an appropriate
sequence). Since $D_k$ is exceptional, we have $Z(t_k\{D_k\})=0$ for
every $k$, but $Z(\alpha)=\alpha$ since $\alpha$ is modified nef by
the above lemma. Consequently, $Z^{-1}(0)$ is not closed, and $Z$ is
not continuous.\\

To get an example of discontinuous Zariski projection, just take $X$ to be the blow-up of $\pet^2$ in at least 9 general points. Such a rational surface is known to carry countably many exceptional curves of the first kind (cf. [Har77], p.409). Since a prime divisor $C$ on a surface is exceptional iff $C^2<0$ (cf. section 4.3), the set of exceptional primes on $X$ is infinite, and we have our example.
 
\subsection{When is a decomposition the Zariski decomposition?}

Suppose that we have a decomposition $\alpha=p+\{N\}$ of a
pseudo-effective class $\alpha$ into the sum of a modified nef class
$p$ and the class of an effective $\er$-divisor $N$. We want a
criterion that tells us when it is the Zariski decomposition of
$\alpha$. We have $N(\alpha)\leq N(p)+N$, and $N(p)=0$ since $p$ is
modified nef, thus $N(\alpha)=N$ happens iff $Z(\alpha)=p$, and our
question is equivalent to the study of the fibers $Z^{-1}(p)$, with $p\in\mnef$.\\   
We will need the following
\begin{defi}[Non-K\"ahler locus] If $\alpha\in\dach$ is a big class,
  we define its non-K\"ahler locus as $E_{nK}(\alpha):=\cap_TE_+(T)$ for
  $T$ ranging among the K\"ahler currents in $\alpha$. 

\end{defi}
Let us explain the terminology:

\begin{theo} Let $\alpha\in\dach$ be a big class. Then:

\noindent(i) The non-nef locus $E_{nn}(\alpha)$ is contained in the
non-K\"ahler locus $E_{nK}(\alpha)$.

\noindent(ii) There exists a K\"ahler current with analytic singularities $T$
in $\alpha$ such that $E_+(T)=E_{nK}(\alpha)$. In particular, the non-K\"ahler locus $E_{nK}(\alpha)$ is an analytic subset of $X$.

\noindent(iii) $\alpha$ is a K\"ahler class iff $E_{nK}(\alpha)$
is empty. More generally, $\alpha$ is a K\"ahler class iff
$\alpha_{|Y}$ is a K\"ahler class for every irreducible component $Y$ of the analytic set $E_{nK}(\alpha)$.
\end{theo}
$Proof$:(i) Since $\alpha$ is big, its non-nef locus $E_{nn}(\alpha)$
is just the set $\{x\in X,\nu(T_{\min},x)>0\}$, since we have
$\nu(\alpha,x)=\nu(T_{\min},x)$ in that case (cf. proposition 3.8). For
every K\"ahler current $T$ in $\alpha$, we have
$\nu(T,x)\geq\nu(T_{\min},x)$ by minimality, and the inclusion
$E_{nn}(\alpha)\subset E_{nK}(\alpha)$ ensues.

(ii) First, we claim that given two K\"ahler currents $T_1$, $T_2$ in
$\alpha$, there exists a K\"ahler current with analytic singularities $T$
such that $E_+(T)\subset E_+(T_1)\cap E_+(T_2)$. Indeed, we can find
$\ep>0$ small enough such that $T_j\geq\ep\omega$. Our currents $T_1$
and $T_2$ thus belong to $\alpha[\ep\omega]$, and admit an infimum
$T_3$ in that set with respect to $\preceq$ (cf. section 2.8). In
particular, $T_3$ is a current in $\alpha$ with $T_3\geq\ep\omega$ and
$\nu(T_3,x)=\min\{\nu(T_1,x),\nu(T_2,x)\}$ for every $x\in X$. By (ii)
of theorem 2.1, there exists a K\"ahler current with analytic
singularities $T$ in $\alpha$ such that $\nu(T,x)\leq\nu(T_3,x)$ for
every $x\in X$, hence $E_+(T)\subset E_+(T_1)\cap E_+(T_2)$, and this
proves the claim.

\noindent Using the claim and (ii) of theorem 2.1, it is easy to
construct a sequence $T_k$ of K\"ahler currents with analytic
singularities such that $E_+(T_k)$ is a decreasing sequence with
$E_{nK}(\alpha)=\cap_kE_+(T_k)$. Since $T_k$ has analytic
singularities, $E_+(T_k)$ is an analytic subset, thus the decreasing
sequence $E_+(T_k)$ has to be stationary (by the strong N\"otherian
property), and we eventually get $E_{nK}(\alpha)=E_+(T_k)$ for some
$k$, as desired. 

(iii) If $\alpha$ is a K\"ahler class, $E_+(\omega)$ is empty for every
K\"ahler form $\omega$ in $\alpha$, and thus so is
$E_{nK}(\alpha)$. Conversely, assume that
$\alpha_{|Y}$ is a K\"ahler class for every component $Y$ of
$E_+(\alpha)$, and let $T$ be a K\"ahler current with analytic
singularities such that $E_+(T)=E_{nK}(\alpha)$. $\alpha$ is then a K\"ahler class by proposition 3.3 of [DP01], qed.\\

We can now state the following 
\begin{theo} Let $p$ be a big and modified nef class. Then the primes $D_1,...,D_r$ contained in the non-K\"ahler locus $E_{nK}(p)$ form an exceptional family $A$, and the fiber of $Z$ above $p$ is the simplicial cone $Z^{-1}(p)=p+V_+(A)$. When $p$ is an arbitrary modified nef class, $Z^{-1}(p)$ is an at most countable union of simplicial cones $p+V_+(A)$, where $A$ is an exceptional family of primes. 
\end{theo}
$Proof$: note that, by the very definitions, for every
pseudo-effective class $\alpha$, the prime components of its negative
part $N(\alpha)$ are exactly the set $A$ of primes $D$ contained in
the non-nef locus $E_{nn}(\alpha)$. Furthermore, $Z(\alpha)+V_+(A)$
is entirely contained in the fiber $Z^{-1}Z(\alpha)$. Indeed, the
restriction of $Z$ to this simplicial cone is a concave map above the
affine constant map $Z(\alpha)$, and both coincide at the relative
interior point $\alpha$, thus they are equal on the whole of
$Z(\alpha)+V_+(A)$. This already proves the last assertion.\\
Assume now that $p$ is modified nef and big, and suppose first that
$\alpha$ lies in $Z^{-1}(p)$. To see that $\alpha$ lies in $p+V_+(A)$,
we have to prove that every prime $D_0$ with $\nu(\alpha,D_0)>0$ lies
in $E_{nK}(p)$, that is: $\nu(T,D_0)>0$ for every K\"ahler current
$T$ in $p$. If not, choose a smooth form $\theta$ in $\{D_0\}$. Since
$T$ is a K\"ahler current, so is $T+\ep\theta$ for $\ep$ small
enough. For $0<\ep<\nu(\alpha,D_0)$ small enough,
$T_{\ep}:=T+\ep\theta+(\nu(\alpha,D_0)-\ep)[D_0]+\sum_{D\neq
  D_0}\nu(\alpha,D)[D]$ is then a positive current in $\alpha$ with
$\nu(T_{\ep},D_0)=\nu(\alpha,D_0)-\ep<\nu(\alpha,D)=\nu(T_{\min},D_0)$
(the last equality holds by proposition 3.8 because $\alpha$ is big
since $p$ is); this is a contradiction which proves the inclusion
$Z^{-1}(p)\subset p+V_+(A)$.\\
In the other direction, let $T$ be a K\"ahler current in $p$, and let
$T=R+\sum\nu(T,D)D$ be its Siu decomposition. $R$ is then a K\"ahler
current with $\nu(R,D)=0$ for every prime $D$, thus its class
$\beta:=\{R\}$ is a modified K\"ahler class. We first claim that we
have $D_j\subset E_{nn}(p-\ep\beta)$ for every $\ep>0$ small
enough and every prime component $D_j$ of the non-K\"ahler locus
$E_{nK}(p)$ of $p$. Indeed, since $p-\ep\beta$ is big for
$\ep>0$ small enough, we have $\nu(p-\ep\beta,D_j)=\nu(T,D_j)$ if $T$
is a positive current with minimal singularities in $p-\ep\beta$, and
we have to see that $\nu(T,D_j)>0$. But $T+\ep R$ is a K\"ahler
current in $p$, thus $D_j\subset E_{nK}(p)\subset E_+(T+\ep R)$ by
definition, which exactly means that $\nu(T+\ep R,D_j)>0$. The claim
follows since $\nu(R,D_j)=0$ by construction of $R$.\\
As a consequence of this claim, each prime $D_1,...,D_r$ of our family
$A$ occurs in the negative part $N(p-\ep\beta)$ for $\ep>0$ small
enough. Consequently, by the first part of the proof, the Zariski
projection of $Z(p-\ep\beta)+\{E\}$ is just $Z(p-\ep\beta)$ for every
effective $\er$-divisor $E$ supported by the $D_j$'s and every $\ep>0$
small enough. Since $p$ is big, $Z$ is continuous at $p$, thus
$Z(p-\ep\beta)$ converges to $Z(p)$, which is just $p$ because the
latter is also modified nef. Finally, $Z$ is also continuous at the
big class $p+\{E\}$, thus the Zariski projection of
$Z(p-\ep\beta)+\{E\}$ converges to that of $p+\{E\}$, and thus
$Z(p+\{E\})=p$ holds. This means that $p+V_+(A)\subset Z^{-1}(p)$, and
concludes the proof of theorem 3.20.   

\subsection{Structure of the pseudo-effective cone}
Using our constructions, we will prove the
\begin{theo}The boundary of the pseudo-effective cone is locally
polyhedral away from the modified nef cone, with extremal rays generated by (the classes of) exceptional prime divisors.
\end{theo}
$Proof$: this is in fact rather straightforward by now: for each prime
$D$, the set $\peff_D:=\{\alpha\in\peff,\nu(\alpha,D)=0\}$ is a closed
convex subcone fo $\peff$. This follows from the fact that
$\alpha\mapsto\nu(\alpha,D)$ is convex, homogeneous, lower
semi-continuous and everywhere non-negative. If
$\alpha\in\partial\peff$ does not belong to $\mnef$, it does not
belong to $\peff_D$ for some prime $D$ by proposition 3.2. For every
$\beta\in\peff$, we have either $\beta\in\peff_D$, or $D$ occurs in
the negative part $N(\beta)$. Therefore, $\peff$ is generated by
$\er_+\{D\}$ and $\peff_D$, and the latter does not contain
$\alpha$. This means that $\partial\peff$ is locally polyhedral near
$\alpha$. Since $\nu(\alpha,D)>0$, we also see that $D$ is exceptional. Finally, the extremal rays of $\peff$ not contained in $\mnef=\cap_D\peff_D$ have to lie outside $\peff_D$ for some exceptional prime $D$, and since $\peff=\peff_D+\er_+\{D\}$, each such extremal ray is generated by $\{D\}$ for some $D$, qed.

\subsection{Volumes}
Recall that the volume of a pseudo-effective class $\alpha$ on a compact K\"ahler $n$-fold is defined to be the supremum $v(\alpha)$ of $\int_XT_{ac}^n$ for $T$ a closed positive $(1,1)$-current in $\alpha$ (cf. [Bou02]). A class $\alpha$ is big iff $v(\alpha)>0$, and the volume is a quantitative measure of its bigness. We have already noticed that $Z(\alpha)$ is big iff $\alpha$ is; we have the following quantitative version:

\begin{prop} Let $\alpha$ be a pseudo-effective class on $X$ compact K\"ahler.
Then $v(Z(\alpha))=v(\alpha)$. 
\end{prop}
The proof is in fact immediate:  if $T$ is a positive current in $\alpha$, then we have $T\geq N(\alpha)$ since $T$ belongs to $\alpha[-\ep\omega]$ for each $\ep>0$, and we deduce that $T\to T-N(\alpha)$ is a bijection between the positive currents in $\alpha$ and those in $Z(\alpha)$. It remains to notice that $(T-N(\alpha))_{ac}=T_{ac}$ to conclude the proof. 

\section{Zariski decomposition on a surface and a hyper-K\"ahler manifold}
It is known since the pioneering work of Zariski [Zar62] that any
effective divisor $D$ on a projective surface admits a unique Zariski
decomposition $D=P+N$, i.e. a decomposition into a sum of
$\ku$-divisors $P$ and $N$ with the following properties:  

(i) $P$ is nef, $N=\sum a_jN_j$ is effective,  

(ii) $P\cdot N=0$,

(iii) The Gram matrix $(N_i\cdot N_j)$ is negative definite.

We want to show that our divisorial Zariski decomposition indeed is a generalization of such a Zariski decomposition on a surface. 

\subsection{Notations}
$X$ will stand for a compact K\"ahler surface, or a compact hyper-K\"ahler manifold. For such an $X$, we denote by $q$ the quadratic form on $\ach$ defined as follows: when $X$ is a surface, we set $q(\alpha):=\int\alpha^2$, and when $X$ is
hyper-K\"ahler, we choose a symplectic holomorphic form $\sigma$, and let
$q(\alpha):=\int\alpha^2(\sigma\overline{\sigma})^{m-1}$ be
the usual Beauville-Bogomolov quadratic form, with $\sigma$ normalized so as to
achieve $q(\alpha)^m=\int_X\alpha^{2m}$ (with $\dim X=n=2m$). In both
cases $(\ach,q)$ is Lorentzian, i.e. it has signature
$(1,h^{1,1}(X)-1)$; the open cone $\{\alpha\in\ach,q(\alpha)>0\}$ has
thus two connected components which are convex cones, and we denote by
$\pos$ the component containing the K\"ahler cone $\ka$. We call
$\pos$ the positive cone (attached to the quadratic form $q$). In
general, given a linear form $\lambda$ on $\ach$, we will denote its
kernel by $\lambda^{\perp}$ and the two open half-spaces it defines by
$\lambda_{>0}$ and $\lambda_{<0}$. The dual $\mathcal{C}^{\star}$ of a
convex cone $\mathcal{C}$ in $\ach$ is seen as a cone in $\ach$, using the duality
induced by $q$.

\subsection{The dual pseudo-effective cone}
In both cases, we shall prove that the modifies nef cone is the dual cone to the pseudo-effective cone.
\subsubsection{The case of a surface}
We suppose that $X$ is a surface. We prove the following essentially well-known
\begin{theo} When $X$ is surface, the K\"ahler cone and the modified K\"ahler cone coincide. The dual pseudo-effective cone is just the nef cone.
\end{theo}
$Proof$: if $\alpha\in\mka$, it can be represented by a K\"ahler current with analytic singularities in codimension 2, that is at some points $x_1,...,x_r$. Therefore we see that the non-K\"ahler locus $E_{nK}(\alpha)$ is a discrete set. Since the restriction of any class to a point is (by convention) a K\"ahler class, theorem 3.19 shows that $\alpha$ lies in fact in $\ka$.\\
Since $\int_X\omega\wedge T$ is positive for every K\"ahler form
$\omega$ and every positive current $T$, we of course have $\ka\subset\cone$, and thus also $\nef=\overline{K}\subset\cone$. The other inclusion is much deeper, since it is a consequence of the Nakai-Moishezon criterion for K\"ahler classes on a surface, as given in [Lam99]. Indeed, this criterion implies that a real $(1,1)$-class $\alpha$ on a K\"ahler surface is a nef class iff $\alpha\cdot\omega\geq 0$ for every $\omega\in\ka$ and $\alpha\cdot C\geq 0$ for every irreducible curve $C$. Since a class in $\cone$ clearly satisfies these conditions, we get $\cone\subset\nef$, and the proof of theorem 4.1 is over.\\

As a consequence, since $\ka$ is contained in $\pos$ and since $\overline{\pos}$ is self dual (just because $q$ is Lorentzian), we get dually that $\overline{\pos}\subset\peff$ and thus that $\pos\subset\peff^0=\bigc$, which means the following: if $\alpha$ is a real $(1,1)$-class with $\alpha^2>0$, then $\alpha$ or $-\alpha$ is big. This generalizes the well known case where $\alpha$ is (the first Chern class of) a line bundle (whose proof is based on Riemann-Roch).

\subsubsection{The hyper-K\"ahler case}
In that case, the dual peudo-effective cone is also equal to the modified nef cone, but the proof uses another description, due to D.Huybrechts, of the dual pseudo-effective cone. In the easy direction, we have:  

\begin{prop}(i) The modified nef cone $\mnef$ is contained in both the dual pseudo-effective cone $\cone$ and the closure of the positive cone $\overline{\pos}$

\noindent(ii) We have $q(D,D')\geq 0$ for any two distinct prime divisors
$D\neq D'$. 
\end{prop}
$Proof$: to prove (i), we only have to prove that
$\mka\subset\cone$. Indeed, $\mka\cap\cone\subset\peff\cap\cone$ is
trivially contained in $\overline{\pos}$. We pick a modified K\"ahler
class $\alpha$ and a pseudo-effective class $\beta\in\peff$, and
choose a K\"ahler current $T$ in $\alpha$ with analytic singularities
in codimension at least 2, and a positive current $S$ in $\beta$. By
section 2.6, the wedge product $T\wedge S$ is well defined as a closed
positive $(2,2)$-current, and lies in the class $\alpha\cdot\beta$. Since $(\sigma\overline{\sigma})^{m-1}$ is a smooth positive form of bidimension $(2,2)$, the integral 
$\int_XT\wedge S\wedge(\sigma\overline{\sigma})^{m-1}$ is positive. But $(\sigma\overline{\sigma})^{m-1}$ is also closed, thus we have 
$$\int_XT\wedge S\wedge(\sigma\overline{\sigma})^{m-1}=\alpha\cdot\beta\cdot\{(\sigma\overline{\sigma})^{m-1}\}=q(\alpha,\beta),$$
so we have proven that $q(\alpha,\beta)\geq 0$ as desired.\\
The second contention is obtained similarly, noting that
$\{D\}\cdot\{D'\}$ contains a closed positive $(2,2)$-current, which
is 
$[D\cdot D']$, where $D\cdot D'$ is the effective intersection cycle.\\

The other direction $\cone\subset\mnef$ is much deeper.
The effective $1$-dimensional cycles $C$ and the effective
divisors $D$ define linear forms on $\ach$ via the intersection form
and the Beauville-Bogomolov form $q$ respectively, and we define a
rational (resp. uniruled) chamber of the positive cone $\pos$ to be a
connected component of $\pos-\cup C^{\perp}$ (resp. $\pos-\cup
D^{\perp}$), where $C$ (resp. $D$) runs over the rational curves
(resp. the uniruled divisors). By a rational curve (resp. a uniruled
divisor) we mean an effective $1$-dimensional cycle all of whose
components are irreducible rational curves (resp. an effective divisor
all of whose components are uniruled prime divisors). The rational
chamber of $\pos$ cut out by all the $C_{>0}$'s (resp. $D_{>0}$)'s will be called the fundamental rational chamber (resp. the fundamental uniruled chamber). When $X$ is a $K3$ surface, the rational and uniruled chambers are the same thing and coincide with the traditional chambers in that situation. We can now state the following fundamental result:
\begin{theo}[Huy99](i) The positive cone $\pos$ is contained in $\peff$.

(ii) If $\alpha\in\pos$ belongs to one of the rational chambers, then there exists a bimeromorphic map $f:X-\to X'$ to a hyper-K\"ahler $X'$ such that 
$$f_{\star}\alpha=\omega'+\{D'\},$$
where $\omega'\in\ka_{X'}$ is a K\"ahler class and $D'$ is an uniruled
$\er$-divisor.\\ 
(iii) When $\alpha\in\pos$ lies in both the fundamental uniruled
chamber and one of the rational chambers, then no uniruled divisor $D'$ occurs in (ii).\\ 
(iv) The fundamental rational chamber coincides with the K\"ahler cone of $X$.
\end{theo}
In fact, [Huy99] states this only for a very general element $\alpha\in\pos$, but we have noticed in [Bou01] that the elements of the rational chambers are already very general in that respect.\\
In the situation (iii), $\alpha$ lies in $f^{\star}\ka_{X'}$ for some
bimeromorphic $f:X-\to X'$ towards a hyper-K\"ahler $X'$. The union of
such open convex cones $\ka_f:=f^{\star}\ka_{X'}$ is called the
bimeromorphic K\"ahler cone, and is denoted by $\bka$. The union in
question yields in fact a partition of $\bka$ into open convex cones
$\ka_f$ (since a bimeromorphic map between minimal manifolds which
sends one K\"ahler class to a K\"ahler class is an isomorphism by a
result of A.Fujiki); $\bka$ is an open cone, but definitely not convex
in general. (iii) tells us that each intersection of a rational
chamber with the fundamental uniruled chamber is contained in $\bka$,
and thus in one of the $\ka_f$'s.\\
We can now describe the dual pseudo-effective cone:
\begin{prop} The dual pseudo-effective $\cone$ of a hyper-K\"ahler manifold coincides with the modified nef cone $\mnef$. 
\end{prop}
$Proof$: by proposition 4.2, it remains to see that $\cone$ is
contained in the modified nef cone $\mnef$. By (i) of theorem 4.3, we
have $\cone\subset\overline{\pos}$, and it will thus be enough to show
that an element of the interior of $\cone$ which belongs to one of the
rational chambers lies in $\mnef$. But an element $\alpha$ of the
interior of $\cone$ has $q(\alpha,D)>0$ for every prime $D$, thus it
certainly lies in the fundamental uniruled chamber. If $\alpha$ lies
in both the interior of $\cone$ and one of the rational chambers, it
therefore lies in $\ka_f=f^{\star}\ka_{X'}$ for some bimeromorphic
$f:X-\to X'$, and it remains to see that $\ka_f\subset\mnef$. But if
$\omega$ is a K\"ahler form on $X'$, its pull-back
$T:=f^{\star}\omega$ can be defined using a resolution of $f$, and it
is easy to check that $T$ is a K\"ahler current with $\nu(T,D)=0$ for
every prime $D$, since $f$ induces an isomorphism $X-A\to X'-A'$ for
$A$, $A'$ analytic subsets of codimension at least 2 (this is because
$X$ and $X'$ are minimal). Therefore, $\{T\}=f^{\star}\{\omega\}$
belongs to $\mka\subset\mnef$, qed.

\subsection{Exceptional divisors}
When $X$ is a surface or a hyper-K\"ahler manifold, the fact that a family $D_1,...,D_r$ of prime divisors is exceptional can be read off its Gram matrix.

\begin{theo} A family $D_1,...,D_r$ of prime divisors is exceptional iff its Gram matrix $(q(D_i,D_j))$ is negative definite.
\end{theo} 
$Proof$: let $V$ (resp. $V_+$) be the real vector space of
$\er$-divisors (resp. effective $\er$-divisors) supported by the
$D_j$'s. We begin with a lemma of quadratic algebra:
\begin{lem} Assume that $(V,q)$ is negative definite. Then every $E\in V$ such that $q(E,D_j)\leq 0$ for all $j$ belongs to $V_+$.
\end{lem}
$Proof$: if $E\in V$ is non-positive against each $D_j$, we write $E=E_+-E_-$ where $E_+$ and $E_-$ are effective with disjoint supports. We have to prove that $E_-=0$, and this is equivalent by assumption to $q(E_-)\geq 0$. But $q(E_-)=q(E_-,E_+)-q(E_-,E)$. The first term is positive because $E_+$ and $E_-$ have disjoint supports, using (ii) of proposition 4.2, whereas the second is positive by assumption on $E$.\\

Let $D_1,...,D_r$ be primes with negative definite Gram matrix. In
particular, we then have that $\{V_+\}\subset\ach$ meets
$\overline{\pos}$ at $0$ only. Since the modified nef cone $\mnef$ is
contained in $\overline{\pos}$ by proposition 4.2, $\{V_+\}$ $a$
$fortiori$ meets the modified nef cone at $0$ only, which means by
definition that $D_1,...,D_r$ is an exceptional family, and this
proves necessity in theorem 4.5.
In the other direction, assume that $D_1,...,D_r$ is an exceptional
family of primes. We first prove that the matrix $(q(D_i,D_j))$ is
semi-negative. If not, we find an $\er$-divisor $E$ in $V$ with
$q(E)>0$. Writing again $E=E_+-E_-$, with $E_+$ and $E_-$ two
effective divisors in $V_+$ with disjoint supports, we have again
$q(E_+,E_-)\geq 0$ by (ii) of proposition 4.2, and thus
$q(E_+)+q(E_-)\geq q(E)>0$. We may therefore assume that $E$ lies in
$V_+$, with $q(E)>0$. But then $E$ or $-E$ is big, and it has to be
$E$ because it is already effective. Its Zariski projection $Z(\{E\})$
is then non-zero since it is also big (by proposition 3.10), and it lies in both $\{V_+\}$ and $\mnef$, a contradiction.\\
To conclude the proof of theorem 4.5, we may assume (by induction)
that the Gram matrix of $D_1,...,D_{r-1}$ is negative definite. If
$(V,q)$ is degenerate, the span $V'$ of $D_1,...,D_{r-1}$ is such that
its orthogonal space $V'^{\perp}$ in $V$ is equal to the null-space of $V$. We then decompose $D_r=E+F$ in the direct sum $V=V'\oplus V'^{\perp}$. Since $q(E,D_j)=q(D_r,D_j)\geq 0$ for $j<r$, lemma 4.6 yields that $E\leq 0$. Therefore, $F=D_r-E$ lies in $V_+$, and is certainly non-zero. We claim that $\{F\}$ is also modified nef, which will yield the expected contradiction. But $F$ lies in the null-space of $V$, and is therefore non-negative against every prime divisor $D$. If $\alpha$ is a pseudo-effective class, we have $q(\{F\},\alpha)=q(\{F\},Z(\alpha))+q(F,N(\alpha))$. The first term is positive since $Z(\alpha)\in\mnef=\cone$, and the the second one is positive because $F$ is positive against every effective divisor. We infer from all this that $\{F\}$ lies in $\cone=\mnef$, and the claim follows.\\

The theorem says in particular that a prime divisor $D$ is negative iff $q(D)<0$. On a K3 surface, an easy and well-known argument using the adjunction formula shows that the prime divisors with negative square are necessarily smooth rational curves with square $-2$. In higher dimension, we have:
\begin{prop} On a hyper-K\"ahler manifold $X$, the exceptional prime divisors are uniruled.
\end{prop}
$Proof$: since $D$ is exceptional, it lies outside $\overline{\pos}=\pos^{\star}$, and we thus find a class $\alpha\in\pos$ lying in one of the rational chambers such that $q(\alpha,D)<0$. By (ii) of theorem 4.3, there exists a bimeromorphic map between hyper-K\"ahler manifolds $f:X-\to X'$ such that $f_{\star}\alpha=\omega'+\sum a_j D_j'$ with $\omega'$ a K\"ahler class, $a_j\geq 0$ and $D_j'$ a uniruled prime divisor. Since the quadratic form is preserved by $f$, we have $0>q(\alpha,D)=q(\omega',f_{\star}D)+\sum a_j q(D_j',f_{\star}D)$, and $q(D_j',f_{\star}D)$ has to be negative for some $j$. But this implies that the two primes $D_j'$ and $f_{\star}D$ coincide, and thus $D=f^{\star}D_j'$ is uniruled since $D_j'$ is. 
  
\subsection{Rationality of the Zariski decomposition}
We want to prove that the divisorial Zariski decomposition is rational
(when $X$ is a surface or a hyper-K\"ahler manifold) in the sense that
$N(\alpha)$ is a rational divisor when $\alpha$ is a rational
class. We first show the following characterization of the divisorial
Zariski decomposition:
\begin{theo} If $\alpha\in\ach$ is a pseudo-effective class, its
  divisorial Zariski decomposition $\alpha=Z(\alpha)+\{N(\alpha)\}$ is
  the unique orthogonal decomposition of $\alpha$ into the sum of a
  modified nef class and the class of an exceptional effective $\er$-divisor.
\end{theo}
$Proof$: we first prove uniqueness: assume that $\alpha=p+\{N\}$ is an orthogonal
decomposition with $p$ a modified nef class and $N$ an effective
exceptional $\er$-divisor. We claim that $N(\alpha)=N$. To see this,
let $D_1,...,D_r$ be the support of $N$; the Gram matrix
$(q(D_i,D_j))$ is negative definite by theorem 4.5, and $p$ is
orthogonal to each $D_j$ because $q(p,N)=0$ and $q(p,D_j)\geq 0$ for
all $j$ since $p$ is a modified nef class. We have $N(\alpha)\leq
N(p)+N$ and $N(p)=0$ since $p$ is modified nef, thus $N(\alpha)\leq
N$. But $N(\alpha)-N$ is supported by primes $D_1,...,D_r$ whose Gram
matrix is negative definite, and
$q(N(\alpha)-N,D_j)=q(p,D_j)-q(Z(\alpha),D_j)$ is non-positive since
$p$ is orthogonal to $D_j$ and $Z(\alpha)$ belongs to
$\mnef=\cone$. Lemma 4.6 thus yields $N(\alpha)\geq N$, and the claim
follows. To prove theorem 4.8, we will show the existence of an
orthogonal decomposition $\alpha=p+\{N\}$ with $p$ a modified nef
class and $N$ an exceptional $\er$-divisor. When this is done, we must
have $N=N(\alpha)$ by the claim, so that that
$\alpha=Z(\alpha)+\{N(\alpha)\}$ is itself an orthogonal decomposition.
\begin{lem} A pseudo-effective class $\alpha$ lies in $\cone$ iff $q(\alpha,D)\geq 0$ for every prime $D$.
\end{lem}
$Proof$: if $\beta$ is a pseudo-effective class, we write $q(\alpha,\beta)=q(\alpha,Z(\beta))+q(\alpha,N(\beta))$. The first term is positive because $Z(\beta)$ lies in $\cone$, and the second one is positive if $q(\alpha,D)\geq 0$ for each prime $D$.

\begin{lem} Let $\alpha$ be a pseudo-effective class and let $D_1,...,D_r$, $E_1,...,E_p$ be two families of primes such that: 

\noindent(i) $q(\alpha,D_j)<0$ and $q(\alpha,E_i)\leq 0$ for every $j$ and $i$.

\noindent(ii) $E_1,...,E_r$ is an exceptional family. 

Then the union of these two families is exceptional. 
\end{lem}
$Proof$: let $F$ be an effective divisor supported by $D_j$'s and $E_i$'s, and assume that $\{F\}$ is a modified nef class. We have to see that $F=0$. But $q(\alpha,F)$ is positive since $F$ is modified nef, thus we see using (i) that $F$ is in fact supported by $E_i$'s, and then (ii) enables us to conclude that $F=0$ as desired.\\ 

At this point, the argument is similar to [Fuj79]. If the
pseudo-effective class $\alpha$ is already in $\cone$, we trivially
have our decomposition. Otherwise, consider the family $A$ of primes
$D$ such that $q(\alpha,D)<0$. That family is exceptional by lemma 4.10
with $E_1,...,E_p$ an empty family, thus $A$ is finite with negative
definite Gram matrix, and is non-empty by lemma 4.9. Let 
$$\alpha=\alpha_1+\{N_1\}$$
be the decomposition in the direct sum $V^{\perp}\oplus V$, where $V\subset\ach$ is spanned by $A$. We claim that $N_1$ is effective and that $\alpha_1$ is pseudo-effective. Since $q(N_1,D)=q(\alpha,D)<0$ for every $D\in A$, lemma 4.6 yields that $N_1$ is effective. We can also write $N(\alpha)=E+F$ where $E$ and $F$ are effective with disjoint supports and $F$ is supported by elements of $A$. Then for every $D\in A$ we have $q(F-N_1,D)\leq q(N(\alpha)-N_1,D)$ since $E$ and $D$ are disjoint, and $q(N(\alpha)-N_1,D)=q(\alpha_1,D)-q(Z(\alpha),D)$ is non-positive because $\alpha_1$ and $D$ are orthogonal and $Z(\alpha)$ lies in $\cone$. We infer from this that $N(\alpha)\geq N_1$ using lemma 4.6, and $\alpha_1=Z(\alpha)+\{N(\alpha)-N_1\}$ is thus pseudo-effective, and this proves our claim.\\
If $\alpha_1$ lies in $\cone$, we have our decomposition by construction; otherwise, we iterate the construction: let $B$ be the non-empty exceptional family of primes $D$ such that $q(\alpha_1,D)<0$. Since $A$ is already exceptional and $q(\alpha_1,D)=0$ for $D\in A$, we infer from lemma 4.10 that the union $A_1$ of $A$ and $B$ is again an exceptional family. We decompose 
$$\alpha_1=\alpha_2+\{N_2\}$$
in the direct sum $V_1^{\perp}\oplus V_1$, where $V_1\subset\ach$ is
spanned by $A_1$. The same arguments as above show in that case also
that $\alpha_2$ is pseudo-effective, and also that $N_2$ is effective
(since $q(N_2,D)=q(\alpha_1,D)\leq 0$ for each $D\in A_1$). But since
$B$ is non-empty, $A_1$ is an exceptional family strictly bigger than
$A$. Since the length of the exceptional families is uniformly bounded
by the Picard number $\rho(X)$ by theorem 3.14, the iteration of the construction has
to stop after $l$ steps, for which we get a class $\alpha_l$ which is
modified nef. The desired decomposition is then obtained by setting
$p:=\alpha_l$ and $N:=N_1+...+N_l$, which is exceptional since it is
supported by elements of $A\cup A_1\cup...\cup A_l=A_l$ (since $A\subset A_1\subset ...\subset A_1$ by construction). This concludes the proof of theorem 4.8.
   
\begin{cor}[Rationality of the Zariski decomposition] The divisorial
  Zariski decomposition is rational in case $X$ is a surface or a
  hyper-K\"ahler manifold. In particular, when $D$ is a
  pseudo-effective divisor on $X$, the modified nef $\er$-divisor
  $P:=D-N(\{D\})$ is rational and such that the canonical inclusion of
  $H^0(X,\oh(kP))$ in $H^0(X,\oh(kD))$ is surjective for every $k$
  such that $kP$ is Cartier. 
\end{cor}
$Proof$: if $\alpha\in NS(X)\otimes\ku$ is a rational class, $N(\alpha)$ is necessarily the image of $\alpha$ by the orthogonal projection $NS(X)\otimes\ku\to V_{\ku}(\alpha)$, where $V_{\ku}(\alpha)$ is the $\ku$-vector space generated by the cohomology classes of the components of $N(\alpha)$. The latter is therefore rational. As to the second part, let $E$ be an element of the linear system $|kD|$. Since the integration current $\frac{1}{k}[E]$ is positive and lies in $\{D\}$, we have $E\geq kN(\{D\})$. But this exactly means that $kN(\{D\})$ is contained in the base scheme of $|kD|$, as was to be shown.

\begin{prop}[Rationality of the volume] If $p\in\ach$ is a modified nef class on $X$, its volume is equal to
$$v(p)=q(p)^m=\int p^{\dim X}.$$
In general, we have $v(\alpha)=\int Z(\alpha)^{\dim X}$; in
particular, the volume of a rational class is rational.  
\end{prop}
$Proof$: we have already proven in proposition 3.22 that
$v(\alpha)=v(Z(\alpha))$, so only the first assertion needs a
proof. We have shown in [Bou02] that the equality $v(p)=\int p^{\dim
  X}$ is always true when $p$ is a nef class, so the contended
equality holds on a surface. In the hyper-K\"ahler case, since we have
chosen the symplectic form $\sigma$ so that $q(\alpha)^m=\alpha^{2m}$
for any class $\alpha$, we just have to prove $v(p)=q(p)^m$ for
$p\in\mnef$. The latter cone is also the closure of the bimeromorphic
K\"ahler cone $\bka$, so we may assume that $p$ lies in
$f^{\star}\ka_{X'}$ for some bimeromorphic map $f:X-\to X'$ between
hyper-K\"ahler manifolds (because both $q$ and the volume are
continuous). But since $f$ is an isomorphism in codimension 1, the
volume is invariant under $f$, and so is the quadratic form $q$, so we
are reduced to the case where $p$ is a K\"ahler class, for which the
equality is always true as we've said above.  

\section{The algebraic approach}
In this section, we would like to show what the constructions we have made become when $\alpha=c_1(L)$ is the first Chern class of a line bundle on a projective complex manifold $X$. The general philosophy is that the divisorial Zariski decomposition of a big line bundle can be defined algebraically in terms of the asymptotic linear series $|kL|$. When $L$ is just pseudo-effective, sections are of course not sufficient, but we are led back to the big case by approximating. For those who are reluctant to assume projectivity too quickly, we remark that a compact K\"ahler manifold carrying a big line bundle is automatically projective.  

\subsection{From sections to currents and back}
Let $L\to X$ be a line bundle over the projective manifold $X$. Each time $L$ has sections
$\sigma_1,...,\sigma_l\in H^0(X,L)$, there is a canonical way to construct a
closed positive current $T\in c_1(L)$ with analytic singularities as follows:
choose some smooth Hermitian metric $h$ on $L$, and consider 
$$\varphi(x):=\frac{1}{2}\log\sum_jh(\sigma_j(x)).$$
Then we define
$T=\Theta_h(L)+dd^c\varphi$, where $\Theta_h(L)$ is the first Chern form of
$h$. One immediately checks that $T$ is positive and independent of the choice
of $h$, and thus depends on the sections $\sigma_j$ only. $T$ has analytic
singularities exactly along the common zero-scheme $A$ of the $\sigma_j$'s, and
its Siu decomposition therefore writes $T=R+D$, where $D$ is the divisor part
of $A$. When $(\sigma_j)$ is a basis of $H^0(X,L)$, we set $T_{|L|}:=T$. Another way to see $T_{|L|}$ is as the pull-back of the Fubiny-Study form on $\pet H^0(X,L)^{\star}=\pet^N$ (the identification is determined by the choice of the basis of $H^0(L)$) by the rational map\\ $\phi_{|L|}:X-\to\pet H^0(X,L)^{\star}$. $T_{|L|}$ is independent of the choice of the basis up to equivalence of singularities, and carries a great deal of information about the linear system $|L|$: the singular scheme $A$ of $T_{|L|}$ is the base scheme $B_{|L|}$ of the linear system $|L|$, the Lelong number $\nu(T_{|L|},x)$ at $x$ is just the so-called multiplicity of the linear system at $x$, which is defined by 
$$\nu(|L|,x):=\min\{\nu(E,x),E\in|L|\}.$$
If a modification $\mu:\ti{X}\to X$ is chosen such that $\mu^{\star}|L|=|M|+F$, where $M$ has non base-point and $F$ is an effective divisor, then $\mu^{\star}T_{|L|}=T_{\mu^{\star}|L|}=T_{|M|}+F$ where $T_{|M|}$ is smooth since $|M|$ is generated by global sections. The so-called moving self-intersection of $L$, which is by definition $L^{[n]}:=M^n$, is thus also equal to $\int_X(T_{|L|})_{ac}^n$.\\  
When $L$ is a big line bundle, we get for each big enough $k>0$ a positive current $T_k:=\frac{1}{k}T_{|kL|}$ in $c_1(L)$. A result of Fujita (cf. [DEL00]) claims that the volume $v(L)$ is the limit of $\frac{1}{k^n}(kL)^{[n]}$, thus we have $v(L)=\lim_{k\to+\infty}\int_XT_{k,ac}^n$.\\
Finally, if $T_{\min}$ is a positive current with minimal singularities in $c_1(L)$, we can choose a singular Hermitian metric $h_{\min}$ on $L$ whose curvature current is $T_{\min}$ (by section 2.4). If $L$ is still big and if for each $k$ we choose the basis of $H^0(kL)$ to be orthonormal with respect to $h_{\min}^{\otimes k}$, then it can be shown that $T_k\to T_{\min}$, and we will see in 5.2 that $\nu(T_k,x)=\frac{1}{k}\nu(|kL|,x)$ converges to $\nu(T_{\min},x)=\nu(c_1(L),x)$. In some sense, the family $T_k$ deriving from $|kL|$ is cofinite $(c_1(L)^+,\preceq)$.\\ 
It should however be stressed that $T_{|kL|}$ will in general $not$ be a K\"ahler current, even if $L$ is big. Indeed, consider the pull-back $L=\mu^{\star}A$ of some ample line bundle $A$ by a blow-up $\mu$. Then $kL$ will be generated by global sections for $k$ big enough, and $T_{|kL|}$ is thus smooth for such a $k$, but not a K\"ahler current, since $L$ is not ample and a smooth K\"ahler current is just a K\"ahler form.\\   

Conversely, to go from currents to sections is the job of the $L^2$ estimates for the $\overline{\partial}$ operator, e.g. in the form of Nadel's vanishing theorem. Recall that the multiplier ideal sheaf $\ih(T)$ of a closed almost positive $(1,1)$-current $T$ is defined locally as follows: write $T=dd^c\varphi$ locally at some $x$. Then the stalk $\ih(T)_x$ is the set of germs of holomorphic functions at $x$ such that $|f|^2e^{-2\varphi}$ is locally integrable at $x$. Then Nadel's vanishing states that if $T$ is a K\"ahler current in the first Chern class $c_1(L)$ of a line bundle $L$, then $H^q(X,\oh(K_X+L)\otimes\ih(T))=0$ for every $q>0$. In particular, if $V(T)$ denotes the scheme $V(\ih(T))$, then the restriction map $$H^0(X,\oh_X(K_X+L))\to H^0(V(T),\oh_{V(T)}(K_X+L))$$ is surjective. This affords a tool to prove generation of jets at some points, using the following lemma (cf. [DEL00]):
\begin{lem}[Skoda's lemma] If $\nu(T,x)<1$, then $\ih(T)_x=\oh_x$. If $\nu(T,x)\geq n+s$, we have $\ih(T)_x\subset\mathcal{M}_x^{s+1}$.
\end{lem}
To illustrate how this works, let us prove the following algebraic characterization of the non-K\"ahler locus:
\begin{prop} If $L$ is a big line bundle, then the non-K\"ahler locus\\ $E_{nK}(c_1(L))$ is the intersection of the non-finite loci $\Sigma_k$ of the rational maps $\phi_{|kL|}$, defined as the union of the reduced base locus $B_{|kL|}$ and the set of $x\in X-B_{|kL|}$ such that the fiber through $x$ $\phi_{|kL|}^{-1}(\phi_{|kL|}(x))$ is positive dimensional somewhere.  
\end{prop}
$Proof$: If $x_1,...,x_r\in X$ lie outside $E_{nK}(c_1(L))$, then we
can find a K\"ahler current $T\in c_1(L)$ with analytic singularities
such that each $x_j$ lies outside the singular locus of $T$. The
latter being closed, there exists a neighbourhood $U_j$ of $x_j$ such
that $\nu(T,z)=0$ for every $z\in U_j$. We artificially force an
isolated pole at each $x_j$ by setting $\ti{T}=T+\sum_{1\leq j\leq
  r}dd^c(\ep\theta_j(z)\log|z-x_j|)$, where $\theta_j$ is a smooth
cut-off function near $x_j$, and $\ep>0$ is so small that $\ti{T}$ is
still K\"ahler. We have $\nu(\ti{T},x_j)=\ep$, whereas $\nu(\ti{T},z)$
is still zero for every $z\neq x_j$ in $U_j$. We now choose a some
smooth form $\tau$ in $c_1(K_X)$, and consider the current
$T_k:=k\ti{T}-\tau$. It lies in the first Chern class of
$L_k:=kL-K_X$, and is certainly still K\"ahler for $k$ big enough. We
also have $\nu(T_k,z)=0$ for every $z\neq x_j$ close to $x_j$, and
$\nu(T_k,x_j)=k\ep$. Given $s_1,...,s_r$, we see that, for $k$ big
enough, each $x_j$ will be isolated in $E_1(T_k)$, whereas
$\ih(T_k)_{x_j}\subset\mathcal{M}_{x_j}^{s_j+1}$, using Skoda's
lemma. Nadel's vanishing then implies that the global sections of $kL$
generate $s_j$-jets at $x_j$ for every $j$. This implies that the
non-finite locus $\Sigma_k$ is contained in $E_{nK}(c_1(L))$.\\
To prove the converse inclusion, we have to find for each $m$ a K\"ahler current $T_m$ in $c_1(L)$ with $E_+(T_m)\subset\Sigma_m$. To do this, we copy the proof of proposition 7.2 in [Dem97]. 
\begin{lem}If $L$ is any line bundle such that the non-finite locus $\Sigma_m$ of $mL$ is distinct from $X$ for some $m$, then, for every line bundle $G$, the base locus of $|kL-G|$ is contained in $\Sigma_m$ for $k$ big enough. 
\end{lem}
We then take $G$ to be ample, and set $T_m:=\frac{1}{k}(T_{|kL-G|}+\omega)$ with $k$ big enough so that $B_{|kL-G|}\subset\Sigma_m$ and $\omega$ a K\"ahler form in $c_1(G)$.\\
To prove lemma 5.3, note that $|mL|$ is not empty, so we can select a modification $\mu:\ti{X}\to X$ such that $\mu^{\star}|mL|=|\ti{L}|+F$, where $|\ti{L}|$ is base point free. It is immediate to check that it is enough to prove the lemma for $\ti{L}$, so we can assume from the beginning that $L$ is base-point free, with $m=1$. We set $\phi:=\phi_{|L|}:X\to\pet^N$ and $\Sigma:=\Sigma_1$. Upon adding a sufficiently ample line bundle to $G$, it is also clear that we may assume $G$ to be very ample. If $x\in X$ lies outside $\Sigma$, the fiber $\phi^{-1}(\phi(x))$ is a finite set, so we can find a divisor $D\in|G|$ which doesn't meet it. Therefore we have $\phi(x)\in\pet^N-\phi(D)$, so that for $k$ big enough there exists $H\in|\oh_{\pet^N}(k)|$ with $H\geq\phi_{\star}D$ which doesn't pass through $\phi(x)$. The effective divisor $\phi^{\star}H-D$ is then an element of $|kL-G|$ which doesn't pass through $x$. The upshot is: for every $x\in X$ outside $\Sigma$, we have $x\in X-B_{|kL-G|}$ for $k$ big enough. By N\"otherian induction, we therefore find $k$ big enough such that $B_{|kL-G|}$ is contained in $\Sigma$, as was to be shown. 

\subsection{Minimal Lelong numbers}
When $L$ is a big $\er$-divisor, we denote by $L_k:=\lfloor kL\rfloor$ the round-down of $kL$, and by $R_k:=kL-L_k$ the fractional part of $kL$. We then consider the sequence $\frac{1}{k}\nu(|L_k|,x)$. It is easily seen to be subadditive, and therefore $\nu(||L||,x):=\lim_{k\to+\infty}\frac{1}{k}\nu(|kL|,x)$ exists. We then prove the following 
\begin{theo} If $L$ is a big $\er$-divisor on $X$ and $\alpha:=\{L\}\in NS(X)_{\er}$, then 
$$\nu(\alpha,x)=\nu(||L||,x)$$
for every $x\in X$.
\end{theo}
$Proof$: let $L=\sum a_jD_j$ be the decomposition of $L$ into its prime components. We choose arbitrary smooth forms $\eta_j$ in $\{D_j\}$, and denote by $\tau_k:=\sum(ka_j-\lfloor ka_j\rfloor)\eta_j$ the corresponding smooth form in $\{R_k\}$. Since $\tau_k$ has bounded coefficients, we can choose a fixed K\"ahler form $\omega$ such that $-\omega\leq\tau_k\leq\omega$ for every $k$. If $E$ is an effective divisor in $|L_k|$, then $1/k([E]+\tau_k)$ is a current in $\alpha[-1/k\omega]$, therefore $\frac{1}{k}\nu(E,x)\geq\nu(T_{\min,1/k},x)$, where $T_{\min,1/k}$ is a current with minimal singularities in $\alpha[-1/k\omega]$, and this yields
$\lim_{k\to\infty}\frac{1}{k}\nu(|L_k|,x)\geq\lim_{k\to\infty}\nu(T_{\min,1/k},x)=\nu(\alpha,x)$. In the other direction, we
use a related argument in [DEL00], Theorem 1.11. The Ohsawa-Takegoshi-Manivel $L^2$
extension theorem says in particular that if we are given a Hermitian line bundle $(A,h_A)$ with sufficiently positive curvature form, then for every pseudo-effective line bundle $G$ and every singular Hermitian metric $h$ on $G$ with positive curvature current $T\in c_1(G)$ and every $x\in X$, the evaluation map
$$H^0(X,\mathcal{O}(G+A)\otimes\mathcal{I}(T))\to\mathcal{O}_x(G+A)\otimes\mathcal{I}(T)_x$$ is surjective, with an $L^2$ estimate independent of $(G,h)$ and $x\in X$.\\ 
We now fix a Hermitian line bundle $(A,h_A)$ with a sufficiently positive curvature form $\omega_A$ to satisfy the Ohsawa-Takegoshi theorem. We select a positive current with minimal singularities $T_{\min}$ in $\alpha$, and also a K\"ahler current $T$ in $\alpha$, which is big by assumption; we can then find almost pluri-subharmonic functions $\varphi_{\min}$ and $\varphi$ on $X$ such that $T_{\min}-dd^c\varphi_{\min}$ and $T-dd^c\varphi$ are smooth. We set $G_k:=L_k-A=kL-R_k-A=(k-k_0)L+(k_0L-R_k-A)$, and fix $k_0$ big enough so that $k_0T-\omega-\omega_A$ is a K\"ahler current. For $k\geq k_0$, the current $T_k:=(k-k_0)T_{\min}+(k_0T-\tau_k-\omega_A)$ is then a positive current in $c_1(G_k)$, thus we can choose for each $k$ a smooth Hermitian metric $h_k$ on $G_k$ such that $T_k$ is the curvature current of the singular Hermitian metric $\exp(-2(k-k_0)\varphi_{\min}-2k_0\varphi)h_k$. Applying the Ohsawa-Takegoshi to $G_k$ equipped with this singular Hermitian metric, we thus get a section $\sigma\in H^0(X,L_k)$ such that 
$$h_k(\sigma(x))\exp(-2(k-k_0)\varphi_{\min}(x)-2\varphi(x))=1$$
and 
$$\int_Xh_k(\sigma)\exp(-2(k-k_0)\varphi_{\min}-2\varphi)dV\leq C_1,$$
where $C_1$ does not depend on $k$ and $x$. If we choose a basis $\sigma_1,...,\sigma_l$ of $H^0(X,L_k)$, we infer from this that
$$\varphi_{\min}(x)+\frac{1}{k-k_0}\varphi(x)=\frac{1}{2(k-k_0)}\log h_k(\sigma(x))$$
$$\leq\frac{1}{2(k-k_0)}\log\sum h_k(\sigma_j(x))+C_2,$$
where $C_2$ does not depend on $x$. The latter inequality comes from the bound on the $L^2$ norm of $\sigma$, since the $L^2$ norm dominates the $L^{\infty}$ norm. Therefore
$$\frac{1}{k-k_0}\nu(|L_k|,x)\leq\nu(\varphi_{\min},x)+\frac{C_3}{k-k_0},$$
where $C_3$ is a bound on the Lelong numbers of $T$. If we let $k\to\infty$ in the last inequality, we get $\nu(||L||,x)\leq\nu(\alpha,x)$ as desired. 

\subsection{Zariski decompositions of a divisor}
The usual setting for the problem of Zariski decompositions is the following:
let $X$ be a projective manifold, and $L$ a divisor on it. One asks when it is
possible to find two $\er$-divisors $P$ and $N$ such that:

(i) $L=P+N$

(ii) $P$ is nef,

(iii) $N$ is effective,

(iv) $H^0(X,kL)=H^0(X,\lfloor kP\rfloor)$ for all $ k>0$, where the
round-down $\lfloor F\rfloor$ of an
$\er$-divisor $F$ is defined coefficient-wise. 

This can of course happen only if $L$ is already pseudo-effective. When this is
possible, one says that $L$ admits a Zariski decomposition (over $\er$ or
$\ku$, depending whether the divisors are real or rational). We want to show
that, for a big divisor $L$, this can be read off the negative part $N(\{L\})$.

\begin{theo} Let $L$ be a big divisor on $X$, and let $N(L):=N(\{L\})$ and $P(L):=L-N(L)$. Then $L=P(L)+N(L)$ is the unique decomposition $L=P+N$ into a modified nef $\er$-divisor $P$ and an effective $\er$-divisor $N$ such that the canonical inclusion $H^0(\lfloor kP\rfloor)\to H^0(kL)$ is an isomorphism for each $k>0$.
\end{theo}
$Proof$: first, we have to check that $H^0(X,kL)=H^0(X,\lfloor kP(L)\rfloor)$. If $E$ is an effective divisor in
the linear system $|kL|$, we have to see that $E\geq\lceil kN(L)\rceil$. But
$\frac{1}{k}[E]$ is a positive current in $\{L\}$, thus $E\geq kN(L)$,
and so $E\geq\lceil kN(L)\rceil$ since $E$ has integer
coefficients.\\
Conversely, assume that $L=P+N$ is a decomposition as in theorem 5.5. We have to show that $N=N(L)$, i.e. $\nu(\{L\},D)=\nu(N,D)$ for every prime $D$. In
view of theorem 5.4, this will be a consequence of the following

\begin{lem} Suppose that a big divisor $L$ writes $L=P+N$, where $P$ is an $\er$-divisor and $N$ is an effective $\er$-divisor such that $H^0(X,kL)=H^0(X,\lfloor kP\rfloor)$ for every $k>0$. Then we have:

(i) If $P$ is nef, then $\nu(||L||,x)=\nu(N,x)$ for every $x\in X$.

(ii) If $P$ is modified nef, then $\nu(||L||,D)=\nu(N,D)$ for every prime $D$.
\end{lem}
$Proof$: the assumption $H^0(X,kL)=H^0(X,\lfloor kP\rfloor)$ means precisely
that for every $E\in|kL|$ we have $E\geq\lceil kN\rceil$, thus
$\nu(|kL|,x)\geq\sum\frac{\lceil ka_j\rceil}{k}\nu(D_j,x)$ if we write $N=\sum
a_jD_j$. We
deduce from this the inequality
$\lim_{k\to\infty}\frac{1}{k}\nu(|kL|,x)\geq\sum a_j\nu(D_j,x)=\nu(N,x)$. To get the converse inequalities, notice that
$$\nu(|kL|,x)\leq\nu(|P_k|,x)+\nu(kN,x)$$
with $P_k:=\lfloor kP\rfloor$ as before; dividing this out by $k$
and letting $k\to+\infty$, we deduce 
$$\lim_{k\to\infty}\frac{1}{k}\nu(|kL|,x)\leq\lim_{k\to\infty}\frac{1}{k}\nu(kN,x)=\nu(N,x)$$ 
when $P$ is nef, since $\nu(\{P\},x)=\lim_{k\to\infty}\frac{1}{k}\nu(P_k,x)$ is then always zero, and similarly with $D$ in place of $x$ when $P$ is modified nef (remark that $P$ is big because $L$ is). This concludes the proof of theorem 5.5.
\begin{cor}[Cutkosky's criterion] Let $L$ be a big divisor on $X$, and
assume that $\nu(\{L\},D)$ is irrational for some irreducible divisor $D$.
Then there cannot exists a modification $\mu:\ti{X}\to X$ such that
$\mu^{\star}L$ admits a Zariski decomposition over $\ku$.
\end{cor} 
$Proof$: if a modification $\mu$ as stated exists, then the negative part
$N(\mu^{\star}L)$ has to be rational by theorem 5.5, and we get a contradiction
using the following easy 
\begin{lem} Let $\alpha$ be a pseudo-effective class on $X$, and let $\mu:\ti{X}\to X$ be a modification. Then we have 
$$N(\alpha)=\mu_{\star}N(\mu^{\star}\alpha).$$
\end{lem}
$Proof$: very easily checked using that a modification is an isomorphism in codimension 1.
\subsubsection{An example of Cutkosky.} 
We propose to analyze in our setting an example due to S.D.Cutkosky [Cut86]
of a big line bundle $L$ on a 3-fold $X$ whose divisorial Zariski
decomposition is not rational, but whose Zariski projection $Z(\{L\})$ is nef. We start from any projective
manifold $Y$ for which $\nef_Y=\peff_Y$. Thus $Y$ might be a
smooth curve or any manifold with nef tangent bundle (cf. [DPS94]). We
pick two very ample divisors $D$ and $H$ on $Y$, and consider
$X:=\pet(\oh(D)\oplus\oh(-H))$, with its canonical projection
$\pi:X\to Y$. If we denote by $L:=\oh(1)$ the canonical relatively
ample line bundle on $X$, then it is well known that 
$$\ach=\pi^{\star}H^{1,1}(Y,\er)\oplus\er L.$$
Since $D$ is ample, $L$
is big, but it won't be nef since $-H$ is not. We are first interested
in the divisorial Zariski decomposition of $L$. We have a hypersurface
$E:=\pet(\oh(-H))\subset X$, and since $D$ has a section, we see that
$E+\pi^{\star}D\in |L|$. Therefore we get $N(L)\leq N(\pi^{\star}D)+E$; but $\pi^{\star}D$ is nef, so has
$N(\pi^{\star}D)=0$, and we deduce $N(L)\leq E$. Consequently,
$N(L)=\mu_LE$ for some $0\leq\mu_L\leq 1$, and $L=Z(L)+\mu_L E$. We
claim that 
$$\mu_L=\min\{t>0, (L-tE)_{|E}\in\nef_E\}.$$
First, we have
$L-tE=\pi^{\star}D+(1-t)E$, and since $\pi^{\star}D$ is nef, we get
that the non-nef locus $E_{nn}(L-tE)$ is contained in $E$ for $0<t<1$. Therefore $L-tE\in\nef_X$ iff
$(L-tE)_{|E}\in\nef_E$. If this is the case, we have $N(L)\leq N(L-tE)+tE=tE$, and thus $t\geq\mu_L$. Conversely, since $L-\mu_L E=Z(L)$ lies in $\mnef$, we get that $Z(L)_{|E}\in\peff_E=\nef_E$ by proposition 2.4 (since $E$ is isomorphic to $Y$ via $\pi$), and we deduce the
equality. Now, notice that the projection $\pi$ induces an isomorphism
$E\to Y$ such that $L$ becomes $-H$ and thus $E_{|E}$ becomes $-D-H$. The
condition
$(L-tE)_{|E}\in\nef_E$ is turned into $-H+t(D+H)\in\nef_Y$, and we get
in the end 
$$\mu_L=\min\{t>0, -H+t(D+H)\in\nef_Y\}$$
The picture can be made more precise:
\begin{prop} (i) The nef cone $\nef_X$ is generated by $\pi^{\star}\nef_Y$ and $L+\pi^{\star}H$. 

(ii) The pseudo-effective cone $\peff_X$ is generated by $\pi^{\star}\nef_Y$ and by $E$. 

(iii) The only exceptional divisor on $X$ is $E$, and the modified K\"ahler cone coincides with the K\"ahler cone. The Zariski projection $Z(\alpha)$ of a pseudo-effective class $\alpha$ is thus the projection of $\alpha$ on $\nef_X$ parallel to $\er_+ E$. 
\end{prop}
$Proof$: given line bundles $L_1,...,L_r$ on a compact K\"ahler manifold $Y$, a class $\alpha=\pi^{\star}\beta$ over $X:=\pet(L_1\oplus...\oplus L_r)$ is nef (resp. pseudo-effective) iff $\beta$ is. A class $\alpha=\oh(1)+\pi^{\star}\beta$ is nef iff $\beta+L_j$ is nef forall $ j$, and $\alpha$ is big iff the convex cone generated by $\beta+L_1,...,\beta+L_r$ meets the big cone of $Y$, which condition is equivalent (by homogeneity) to: $\beta+$conv$(L_1,...,L_r)$ meets the big cone; finally $\alpha$ is pseudo-effective iff $\beta+$conv$(L_1,...,L_r)$ meets $\peff_Y$. In our case $\alpha=\pi^{\star}\beta+L$ is thus nef iff $\beta-H$ is nef, and $\alpha$ is pseudo-effective iff $\alpha+[-H,D]$ meets $\nef_Y$. The latter condition is clearly equivalent to $\alpha-D\in\nef_Y$. Now an arbitrary class $\alpha$ on $X$ uniquely writes $\alpha=tL+\pi^{\star}\beta$. If $\alpha$ is pseudo-effective, then $t\geq 0$ (since $L$ is relatively ample); if $t=0$, then $\alpha\in\pi^{\star}\nef_Y$. Otherwise, we may assume by homogeneity that $t=1$, and thus (i) and (ii) follow from the above discussion.\\
By (ii), a pseudo-effective class $\alpha$ writes
$\pi^{\star}\beta+tE$ with $\beta$ nef. Therefore we get
$N(\alpha)\leq tE$, and $E$ is thus the only exceptional divisor on
$X$. In fact, we even have $E_{nn}(\alpha)\subset E$, and thus $\alpha$ is nef iff $\alpha_{|E}$ is nef. In particular, we see that $\mka=\ka$ as desired (use proposition 2.4 again).\\

We now assume that $Y$ is a surface. The assumption $\nef_Y=\peff_Y$
implies that $\nef_Y=\overline{\pos}_Y=\peff_Y$, and $\mu_L$ is none but the least of the two roots of the quadratic polynomial in $t$ $(-H+t(D+H))^2$; it will thus be irrational for most choices of $H$ and $D$ (on, say, an
abelian surface). This already yields that the divisorial Zariski decomposition of the rational class $c_1(L)$ will not be rational in general, that is, the analogue of corollary 4.11 is not true in general on a 3-fold.\\
Since $Z(L)$ is nef, the volume of $L$ is just $v(Z(L))=Z(L)^3$, with $Z(L)=(1-\mu_L)L+\mu_L\pi^{\star}D$. The cubic intersection form is explicit on $\ach$ from the relations 
$$L^3-\pi^{\star}(D-H)\cdot L^2-D\cdot H\cdot L=0$$
and $\pi_{\star}L=1$, $\pi_{\star}L^2=D-H$, thus we can check that $v(L)$ is an explicit polynomial of degree 3 in $\mu_L$ which is also irrational for most choices of $D$ and $H$. We conclude: there exists a big line bundle on a projective 3-fold with an irrational volume, by contrast with proposition 4.12.

\section{References.}
\begin{itemize}
\item{\bf [Bou01]} Boucksom, S.\ --- {\sl Le c\^one k\"ahlerien d'une vari\'et\'e hyperk\"ahlerienne}, C.R.A.S. (2001), --.
\item{\bf [Bou02]} Boucksom, S.\ --- {\sl On the volume of a line
    bundle}, math.AG/0201031 (2002).
\item{\bf [Cut86]} Cutkosky, S.D.\ --- {\sl Zariski decomposition of divisors on algebraic varieties}, Duke Math. J. {\bf 53} (1986), 149--156.
\item{\bf [Dem82]} Demailly, J.-P.\ --- {\sl Estimations $L^2$ pour l'op\'erateur $\overline{\partial}$ d'un fibr\'e vectoriel holomorphe semi-positif au dessus d'une vari\'et\'e k\"ahlerienne compl\`ete}, Ann. Sci. Ecole Norm. Sup. {\bf 15} (1982), 457--511.
\item{\bf [Dem92]} Demailly, J.-P.\ --- {\sl Regularization of closed positive
currents and intersection theory}, J. Alg. Geom. {\bf 1} (1992), 361--409.
\item{\bf [Dem97]} Demailly, J.-P.\ --- {\sl Algebraic criteria for Kobayashi hyperbolic projective varieties and jet differentials}, Proc. Symp. Pure Math. {\bf 62.2} (1997).
\item{\bf [DPS94]} Demailly, J.-P.; Peternell, T.; Schneider, M.\ ---
  {\sl Compact complex manifolds with numerically effective tangent bundles},
J. Alg. Geom. {\bf 3} (1994), 295-345.
\item{\bf [DEL00]} Demailly, J.-P.; Ein, L.; Lazarsfeld, R.\ --- {\sl A
subadditivity property of multiplier ideals}, math.AG/0002035 (2000). 
\item{\bf [DPS00]} Demailly, J.-P.; Peternell, T.; Schneider, M.\ ---
  {\sl Pseudoeffective line bundles on compact K\"ahler manifolds},
math.AG/0006205 (2000). 
\item{\bf [Fuj79]} Fujita, T.\ --- {\sl On Zariski problem},
  Proc. Japan Acad., Ser. A {\bf 55} (1979), 106--110.
\item{\bf [Fuj89]} Fujita, T.\ --- {\sl Remarks on quasi-polarized varieties},
Nagoya Math. J. {\bf 115} (1989), 105--123.
\item{\bf [Har77]} Hartshorne, R.\ --- {\sl Algebraic geometry}, Springer Verlag, GTM {\bf 52} (1977).
\item{\bf [Huy99]} Huybrechts, D.\ --- {\sl The K\"ahler cone of a
    compact hyperk\"ahler manifold}, math.AG/9909109 (1999).
\item{\bf [Lam99]} Lamari, A.\ --- {\sl Courants k\"ahleriens et surfaces
compactes}, Ann. Inst. Fourier {\bf 49} (1999), 249--263.
\item{\bf [Pau98]} Paun, M.\ --- {\sl Sur l'effectivit\'e num\'erique
    des images inverses de fibr\'es en droites}, Math. Ann. {\bf 310}
  (1998), 411--421.
\item{\bf [Siu74]} Siu, Y.T.\ --- {\sl Analyticity of sets associated to Lelong numbers and the extension of closed positive currents}, Invent. Math. {\bf 27}(1974), 53--156.
\item{\bf [Zar62]} Zariski, O.\ --- {The theorem of Riemann-Roch for high multiples of an effective divisor on an algebraic surface}, Ann. of Math. (2) {\bf 76} (1962), 560--615.
\end{itemize}

\end{document}